\documentclass[12pt]{article}

\usepackage{mathrsfs}
\usepackage{amsmath}
\usepackage{amsfonts}
\usepackage{amssymb}
\usepackage[all,cmtip]{xy}
\usepackage{hyperref}
\usepackage{enumerate}  
\usepackage{ulem}
\usepackage{tikz-cd}
\usepackage[all,cmtip]{xy}
\usepackage[toc,page]{appendix}
\usepackage{amsmath, amsthm, amsfonts, amssymb}
\usepackage{graphicx}
\usepackage{subcaption}

\def\<{\langle}
\def\>{\rangle}

\def\-{\overline}

\def\endpf{\hbox{\vrule height1.5ex width.5em}}

\def\endpf{\hbox{\vrule height1.5ex width.5em}}

\def\-{\overline}

\def\endpf{\hbox{\vrule height1.5ex width.5em}}

\def\endpf{\hbox{\vrule height1.5ex width.5em}}

\def\-{\overline}

\newtheorem{question}{Question}[section]
\newtheorem{theorem}{Theorem}[section]
\newtheorem{lemma}[theorem]{Lemma}
\newtheorem{corollary}[theorem]{Corollary}
\newtheorem{proposition}[theorem]{Proposition}

\newtheorem*{theorem*}{Theorem}
\date{}
\theoremstyle{definition}
\newtheorem{definition}[theorem]{Definition}
\newtheorem{remark}[theorem]{Remark}
\newtheorem{example}[theorem]{Example}
\begin{document}

\title{\bf A geometric criterion for prescribing residues and some applications}

\author{Hanlong Fang}

\maketitle
\begin{abstract}
Let $X$ be a compact complex manifold and $D$ be a $\mathbb C$-linear form sum of divisors of $X$. A theorem of Weil and Kodaira says that  if $X$ is K\"ahler, then there is a closed logarithmic $1$-form with residue divisor $D$  if and only if $D$ is homologous to zero in $H_{2n-2}(X,\mathbb C)$. We generalized their theorem to general compact complex manifolds. The necessary and sufficient condition is described by a new invariant called $\mathcal Q$-flat class. In the second part of the paper, we classify all the pluriharmonic functions on a compact algebraic manifold with mild singularities.
\end{abstract}

\tableofcontents

\bigskip\bigskip
\section{Introduction}
In this paper, we study the following two questions for compact complex manifolds.
\begin{question}[Inverse residue problem]\label{Q1}
Find closed meromorphic $1$-forms (called Abelian diffential in dimension one) with given residues.
\end{question}
\begin{question}[Existence of pluriharmonic functions]\label{Q2} Construct and classify pluriharmonic functions  locally taking the form of
	\begin{equation}\label{log}
	g_1(z)+g_2(\bar z)+\sum_{i=1}^la_i\log |f_i|^2,
	\end{equation}
	where $a_1,\cdots,a_l\in\mathbb C$ and $g_1,g_2,f_1,\cdots,f_l$ are meromorphic functions.
\end{question}

Integrals of meromorphic $1$-forms on Riemann surfaces played an important role in the development of the theory of complex analysis in one variable. There are two well-known holomorphic invariants associated with each meromorphic $1$-form called residues and pole orders as follows. Suppose a meromorphic $1$-form $g$ has the Laurent expansion near point $\xi=0$ in a local chart $(U,\xi)$ as 
\begin{equation}
g=h(\xi)d\xi=\big(\sum\limits_{i=-l}^{\infty}c_i\xi^i\big)d\xi\,\,,\,\,\,c_{-l}\neq 0.
\end{equation}
Then $l$ is the pole order and   $c_{-1}$ is the residue at $\xi$.  This leads to the following classification of meromorphic $1$-forms by types (see \cite{S2} for details):
\begin{enumerate}
	\item A differential is called an {\bf Abelian differential of  the first kind} if it is regular on the Riemann surface, that is, if it has no poles.
	\item A differential is called an {\bf Abelian differential of  the second kind} if it has at least one pole and if, in addition, its residue at each pole is zero.
	\item A differential is called  an {\bf Abelian differential of  the third kind} if it has at least one nonzero residue.
\end{enumerate}

Regarding Question \ref{Q1}, the following classical theorems give a thorough understanding of the compact Riemann surfaces case.

\begin{theorem}[Theorem 1 in \S 4.1 of \cite{S2}]\label{res}
Let $X$ be a compact Riemann surface.  The sum of the residues of an Abelian differential over $X$ is always zero.
\end{theorem}

\begin{theorem}[Theorem 1 in \S 4.2 of \cite{S2}]\label{e1} For each point $p$ on a Riemann surface $X$ and $l=2,3,\cdots$ there exist Abelian differentials of the second kind $dE_l(p)$ and $dF_l(p)$ with a single pole at $p$ with pole order $l$. Furthermore, all the periods of the integral $E_l(p)$ are pure imaginary and those of $F_l(p)$ are real.
\end{theorem}

\begin{theorem}[Theorem 2 in \S 4.2 of \cite{S2}]\label{e2} For any two distinct points $p$ and $q$ on a Riemann surface $X$ there exists an Abelian differential $dE(p,q)$  of the third kind which is regular apart from simple poles $p$ and $q$ with residues $1$ and $-1$. Also, all the periods of the integral $E(p,q)$ are pure imaginary.
\end{theorem}

Picard and Lefschetz generalized the concept of Abelian differentials on algebraic surfaces.   Notice that an Abelian differential on Riemann surfaces is closed, and hence the line integrals are homotopic invariant and the residues are well-defined. When the complex  dimension of a complex manifold is larger than or equal to two, the closeness property does not hold automatically. Therefore, in order to attach residues to a differential, it is natural to make the closeness property as an additional  assumption.  

In higher dimension, Hodge and Atiyah \cite{HA} generalized the concept of Abelian differentials by sheaf theory (see \S 2 for more details).   Let $X$ be a complex manifold and $W$ be a reduced divisor $W$ of $X$, namely, an anlytic subvariety of $X$ of codimension one.  Denote by $H^0(X,R^1(W))$ the space of the residue divisors supported on $W$ (an analogue of residues), or equivalently,  $\mathbb C$-linear formal sums of divisors of $X$ supported on $W$; denote by $H^0(X,\Phi^1(*W))$ the space of the closed meromorphic $1$-forms with poles on $W$; denote by $H^0(X,d\Omega^0(*W))$ the subspace of $H^0(X,\Phi^1(*W))$ consisting of the closed meromorphic $1$-forms  with trivial residue divisors (an analogue of Abelian differentials of the second kind). Then, they derived the following long exact sequences:
\begin{equation}\label{iH1d}
\begin{split}
\rightarrow H^0(X,&\Omega^0(*W))\rightarrow H^0(X,d\Omega^0(*W))\rightarrow H^1(X,\mathbb C)\rightarrow H^1(X,\Omega^0(*W))\rightarrow \\
&\rightarrow H^1(X,d\Omega^0(*W))\rightarrow H^2(X,\mathbb C)\rightarrow H^2(X,\Omega^0 (*W))\rightarrow\cdots,
\end{split}
\end{equation}
\begin{equation}\label{iH0R}
0\rightarrow H^0(X,d\Omega^0(*W))\rightarrow H^0(X,\Phi^1(*W))\xrightarrow{\rm Res} H^0(X,R^1(W))\xrightarrow{\Delta^0} H^1(X,d\Omega^0(*W))\rightarrow\cdots.
\end{equation}

As a consequence of the long exact sequences $(\ref{iH1d})$ and $(\ref{iH0R})$, one derives immediately the following abstract criterion for Question \ref{Q1}.
\begin{theorem}[\cite{HA}]\label{HA}Suppose $X$ is a compact complex manifold and $W$ is a reduced divisor of $X$. Let $D$ be an element of $H^0(X,R^1(W))$. Then there is a closed meromorphic $1$-form with residue divisor $D$ if and only if $\Delta^0(D)=0$ in $H^1(X,d\Omega^0(*W))$.
\end{theorem}

On the other hand, when $X$ is a compact K\"ahler manifold, Weil and Kodaira  gave a geometric criterion for the existence of a closed  logarithmic $1$-form  with prescribed residues by using harmonic integrals and potential theory. (See \S 4.2 for the definition of logarithmic forms.) We rewrite their theorem as follows.
\begin{theorem}[\cite{W} and \cite{Ko}\label{WK}]Let $X$ be a compact K\"ahler manifold of complex dimension $n$ and $W$ be a reduced divisor of $X$. Let $D$ be an element of $H^0(X,R^1(W))$. Then there is a closed logarithmic $1$-form with residue divisor $D$ and with poles on $W$  if and only if $D$ is homologous to zero in $H_{2n-2}(X,\mathbb C)$.
\end{theorem}

In the first part of this paper, we investigate Question \ref{Q1} following Hodge and Atiyah's sheaf theoretical method and generalize the geometric criterion of Weil and Kodaira to general compact complex manifolds. 

Before proceeding, first recall the following fact (see Proposition 3.1 in \cite{P} for instance).  Let  $D$ be a reduced divisor of a complex manifold $X$. Suppose   $D$ is a flat divisor in the sense that  under a certain trivialization the transition functions of $[D]$, which is the line bundle associated with $D$,  can be taken as constant functions. Then, there exists a closed meromorphic $1$-form
with simple poles along the support $D$ and holomorphic on the complement.

We show that the geometric criterion for general compact complex manifolds is given a holomorphic invariant measuring the flatness of divisors as follows.  (See \cite{F} for the definition of the $\mathcal Q$-flat class of a holomorphic line bundle; see \S 3.4 for the definition of the $\mathcal Q$-flat class of a $\mathbb C$-linear formal sum of divisors.)

\begin{theorem}\label{ri}Let $X$ be a compact complex manifold and $W$ be a reduced divisor of $X$. Let $D\in H^0(X,R^1(W))$. Then the following statements are equivalent:
	\begin{enumerate}[]
		\item $\bullet$  the $\mathcal Q$-flat class of $D$ is trivial in $H^1(X,d\Omega^0)$;
		\item $\bullet$  there is a closed logarithmic $1$-form with residue divisor $D$ and with poles on $W$. 
	\end{enumerate}
\end{theorem}
Since the $\mathcal Q$-flat class of a holomorphic line bundle is trivial if only if the line bundle is flat up to some positive multiple (see Theorem 1.11 in \cite{F}),  we have the following interesting corollary of Theorem \ref{ri}.
\begin{corollary}\label{intc}
Let $X$ be a compact complex manifold and $W$ be a reduced divisor of $X$.  If there exist some positive integer $m$ such that $mW$ is flat, then there exists a closed meromorphic $1$-form with simple poles along the support $W$ and holomorphic on the complement.
\end{corollary} 

Next, we study the closed meromorphic $1$-forms with poles of arbitrary order by refining Hodge and Atiyah's criterion (see Theorem \ref{knc}). In particular, we derive the following topological constraint on the residue divisors.
\begin{theorem}\label{nc} Let $X$ be a compact complex manifold. Then the residue divisor of a closed meromorphic $1$-form on $X$  is homologous to zero in $H_{2n-2}(X,\mathbb C)$.
\end{theorem}

We say that a compact complex manifold $X$ has Property $(H)$ (see Definition \ref{ph}) if and only if
\begin{equation}
\dim H^1(X,\mathbb C)=\dim H^{0}(X,d\Omega^0)+\dim H^1(X,\mathcal O_X).
\end{equation} 
In the following, we apply the above general results to compact complex manifolds with Property $(H)$.  We first reduce the holomorphic criterion in Theorem \ref{ri} to a topological criterion, which is similar to that in Theorem \ref{WK}.
\begin{theorem}\label{imain}  Assume that $X$ is a compact complex manifold with Property $(H)$ and that $W$ is a reduced divisor of $X$. Let $D\in H^0(X,R^1(W))$.  Then the following statements are equivalent:
	\begin{enumerate}[]
		\item $\bullet$  $D$ is homologous to zero in $H_{2n-2}(X,\mathbb C)$;
		\item $\bullet$  there is a closed logarithmic $1$-form with residue divisor $D$ and with poles on $W$.
	\end{enumerate}
\end{theorem}

In particular, Theorem \ref{imain} gives an alternate proof of Theorem \ref{WK}, for compact K\"ahler manifolds have Property $(H)$.

Recall a well known result of Deligne \cite{D} that each logarithmic form  on projective manifolds is closed. This result later was generalized to K\"ahler manifolds and  complex manifolds of Fujiki class $\mathcal C$ by  Noguchi \cite{No} and Winkelmann \cite{Wi}, respectively.  We further generalize it for $1$-forms in the following way. 
\begin{theorem}\label{de}
	Let $X$ be a compact complex manifold with Property $(H)$. Let $W$ be an effective reduced divisor on $X$. Then each logarithmic $1$-form on $X$ is a sum of a holomorphic $1$-form and a closed logarithmic $1$-form, that is, 
	\begin{equation}
	H^0(X,\Omega^1(\log W))=H^0(X,\Omega^1 )+H^0(X,\Phi^1(W)).
	\end{equation}
\end{theorem}

We also have the following decompostion result for closed meromorphic $1$-forms with poles of arbitrary order.
\begin{theorem}\label{df}
	Let $X$ be a compact complex manifold with Property ($H$). Then every closed meromorphic $1$-form is a sum of a closed logarithmic $1$-form and a closed meromorphic $1$-form of the second kind.
\end{theorem}
\medskip

In the second part of the paper,  we turn to the study of pluriharmonic functions on projective manifolds (Question \ref{Q2}). Recall that a pluriharmonic function $f$ (possibly singular) is a solution of the following overdetermined system of partial differential equations:
\begin{equation}\label{plur}
\partial\overline \partial {f}=0.
\end{equation}
The only regular solutions of equation $(\ref{plur})$  on a compact complex manifold are constant functions; meromorphic functions and anti-meromorphic functions are its singular solutions. We first derive the following theorem.
\begin{theorem}\label{iepluri}
	
	Let X be a compact algebraic manifold. Assume that $W$ is a reduced divisor of $X$ and   there is an effective, ample divisor of $X$ whose support is contained in $W$.  Then for every closed meromorphic $1$-form with poles on $W$, there exists a closed anti-meromorphic $1$-form on $X$ with poles on $W$, so that the  integral of the sum of these two differentials is a single-valued function on $X\backslash W$.  In particular,  the integral is a pluriharmonic function with singularities on $W.$ 
\end{theorem}

As an application of Theorem \ref{iepluri}, we have the following theorem classifying all the singular solutions of equation $(\ref{plur})$ with local form $(\ref{log})$.
\begin{theorem}\label{121} Let X be a compact algebraic manifold. Denote by $K(X)$ the vector space of meromorphic functions on $X$; denote by $\overline K(X)$ the vector space of anti-meromorphic functions on $X$; denote by $dK(X)$ the vector space of the differentials of meromorphic functions on $X$.  Denote by $Ph(X)$ the vector space of the pluriharmonic functions on $X$ of local form $(\ref{log})$; denote by $Ph_0(X)$ the vector space of the pluriharmonic functions on $X$ of local form $(\ref{log})$ without log terms. Then the following natural homomorphisms induced by differentiation are isomorphisms:
\begin{equation}
\kappa:Ph(X)/(K(X)+\overline K(X))\xrightarrow{\cong} H^0(X,\Phi^1(*))/dK(X)\,\,;
\end{equation}	
\begin{equation}
\kappa_0:Ph_0(X)/(K(X)+\overline K(X))\xrightarrow{\cong} H^0(X,d\Omega^0(*))/dK(X)\,\,.
\end{equation}
\end{theorem}

\bigskip
We now briefly describe the organization of the paper and the basic ideas for the proof of theorems. A nature approach to prove Theorem $\ref{ri}$ and Theorem $\ref{imain}$ is to interpret sheaf cohomology  as \v Cech cohomology. In order to establish the isomorphism between these two cohomologies of $X$, we shall prove that  certain  cohomology groups are trivial. However, the Hodge-Atiyah exact sequence is not a sequence of coherent $\mathcal O_X$-sheaves, and hence Cartan theorem B does not apply. This difficulty is settled by three lemmas: truncation lemma, good cover lemma and acyclic lemma. Next, we construct  explicitly the double delta map from the residue divisor group  to the obstruction group by diagram chasing and prove Theorem $\ref{ri}$ and Theorem $\ref{imain}$. Then, we show that, for manifolds with Property $(H)$, the  criterion for the existence of logarithmic $1$-forms with given residues coincides with the criterion for the  existence of closed meromorphic $1$-forms with given residues. Therfore,  each logarithmic $1$-form can be decomposed into two parts as in Theorem \ref{de}. Similar to the logarithmic case, we derive the criterions for prescribing residues for closed meromorphic $1$-forms with arbitrary pole order and as a consequence of which we prove Theorem \ref{nc}. 

In the second part of the paper, we investigate Question \ref{Q2}. Firstly, we describe the two kinds of obstructions  for getting a single-valued function by integrating closed meromorphic $1$-forms, namely, the long period vectors corresponding to the integrals along loops in $H_1(X,\mathbb C)$, and the short period vectors corresponding to  the integrals along small loops around irreducible components of the residue divisor.  Since such an integral is single-valued if and only if all periods vanish and it is impossible to carry out a cancellation of periods merely in the holomorphic category except the trivial case,  we produce closed anti-meromorphic $1$-forms with opposite periods, and then sum up the pairs to get a single-valued function.  This procedure is made explicitly in the language of gardens and pairs{\footnote{\label{myfootnote}These are named after classical gardens of Suzhou.}}. At last,  we will prove Theorem \ref{iepluri} and Theorem $\ref{121}$ by a careful cancellation of the periods. 

The organization of the paper is as follows: In \S 2, we introduce Hodge and Atiyah's sheaf theoretical method. In \S 3.1, we  establish a special truncation of Hodge-Atiyah exact sequence. In \S 3.2, we prove the acyclic lemma. In \S 3.3, diagram chasing method is used to show the geometric meaning of the double delta map. In \S 3.4, we define the $\mathcal Q$-flat class of a $\mathbb C$ linear formal sums of divisors and derive some basic properties of manifolds with the Property $(H)$.   In \S 4.1, we prove Theorem \ref{ri} and Theorem \ref{imain}. In \S 4.2, we prove Theorem \ref{de}.
In \S 5, we prove Theorem \ref{nc}. In \S 6, we prove Theorem $\ref{iepluri}$ and Theorem \ref{121}. 

For reader's convenience, we include in Appendix I a detailed proof for the existence of a (very) good cover for a compact complex manifold which is crucial in the comparison of two cohomologies. Also, we include in Appendix II a proof for the existence of a smooth, transversal two-chain which is used for calculating periods.
\medskip

\section*{Acknowledgement} The author  appreciates greatly his advisor Prof. Xiaojun Huang for the inspiring course on Riemann surfaces. He thanks Zhan Li for his reading of the draft and making suggestions with patience. Also, he would like to thank Xu Yang who provides a physical explanation of the logarithmic poles and thanks Song-Yan Xie for his help with Latex and English.  Finally, he thanks Xin Fu, Qingchun Ji and Jian Song for bringing up to him the topic of constructing plurharmonic functions.

\section{Preliminary and the Hodge-Atiyah sequences}

Let $X$ be a complex manifold and $\Omega^1_X$ (or $\Omega^1$ for short) the cotangent bundle of $X$. For each meromorphic $1$-form $f$ on $X$, we define an ideal sheaf on $X$ characterizing the singularities of $f$ as follow.

\begin{definition} Define a presheaf $\mathcal P_f^{pre}$ as follows: for each open set (in the Euclidean topology) $U\subset X$,
\begin{equation}
\Gamma(\mathcal P_f^{pre},U):=\{h\in\mathcal O_X(U)\big|\,\,hf|_{U}\in H^0(U,\Omega^1)\}.
\end{equation}
The sheafification of this presheaf is called the denominator ideal sheaf associated with $f$ and denoted  by $\mathcal P_f$.	
	
\end{definition}

\begin{lemma}\label{div}The denominator ideal sheaf associated with $f$ is locally free and of rank 1; that is, $\mathcal P_f$ defines a divisor.
	
\end{lemma}
\noindent{\bf{Proof :}} We choose an open set $U$ of $X$ with complex coordinates $(z_1,\cdots,z_n)$ such that
\begin{equation}
f=\frac{f_1}{g_1}dz_1+\frac{f_2}{g_2}dz_2+\frac{f_3}{g_3}dz_3+\cdots+\frac{f_n}{g_n}dz_n,
\end{equation}
where $f_i,g_i$ are holomorphic functions over $U.$ Without loss of generality, we  assume  $f_i$ and $g_i$ are coprime for  $i=1,\cdots,n.$ Let $h$ be the least common multiple of $g_1,g_2,\cdots,g_n.$ It is easy to verify that $\mathcal P_f(U)=(h).$ We complete the proof of the lemma. $\endpf$ 
\medskip

Following \cite{HA}, we introduce some notations and recall some important results therein. 
In the remaining of this section we assume that $W$ is a reduced  divisor on $X$, namely, $W$ is an analytic subvariety of $X$ of complex codimension one; we also assume that all the open sets are in the Euclidean topology.  
  
For each integer $q\geq 0$, denote by $\Omega^q(kW)$ the sheaf of germs of meromorphic $q$-forms having, as their only singularities, poles of order at most $k$ on the components of $W$. (We view meromorphic functions as $0$-forms, when $q=0$.)   The union  of the sheaves $\Omega^q(kW)$ as $k\rightarrow\infty$ we denote by $\Omega^q(*W)$; it is just the sheaf of germs of meromorphic $q$-forms with poles of any order on $W$. Similarly the union of the sheaves $\Omega^q(*W)$ as $W$ runs through all  reduced  divisors of $X$ is denoted by $\Omega^q(*).$

Define a presheaf  by
\begin{equation}
d\Omega^{q}(kW)(U):=\{df|f\in \Omega^{q}(kW)(U)\}\,\,\text{for each open set\,\,} U;
\end{equation}  
the $\mathbb C$-sheaf $d\Omega^{q}(kW)$ is the sheafification of this presheaf.  Denote by $d\Omega^{q}(*W)$ the union of of the sheaves $d\Omega^{q}(kW)$ as $k\rightarrow\infty$; denote by $d\Omega^{q}(*)$ the union of the sheaves $d\Omega^{q}(*W)$ as $W$ runs through all  reduced divisors of $X$.  Denote  by $\Phi^q(kW)$, $\Phi^q(*W)$ and $\Phi^q(*)$   the subsheaves of $\Omega^q(kW)$, $\Omega^q(*W)$ and $\Omega^q(*)$, respectively, consisting of germs of closed forms.   Moreover, define the sheaf $R^q(W)$ and $R(*)$ by the following exact sequences, respectively:
\begin{equation}\label{*W}
\begin{split}
0\rightarrow &d\Omega^{q-1}(*W)\rightarrow\Phi^q(*W)\rightarrow R^q(W)\rightarrow 0\,;\\
0&\rightarrow d\Omega^{q-1}(*)\rightarrow\Phi^q(*)\rightarrow R^q(*)\rightarrow 0\,.
\end{split}
\end{equation}

Next we will define $\mathbb C$-sheaves $\mathcal D(W)$ and $\mathcal D(*)$.  Let $U$ be an open set of $X$. Denote by $\{W^U_{h}\}$ the irreducible components of $W$ in $U$ and by ${\mathbb C}_{W^U_h}$ the constant sheaf on $W^U_h$. Since ${\mathbb C}_{W^U_h}$  can be viewed as a $\mathbb C$-sheaf on $U$, we define a presheaf ${\mathcal D^{pre}(W)}$ by 
\begin{equation}
\Gamma({\mathcal D^{pre}(W)},U):=\sum_h{\mathbb C}_{W^U_h}=\{\sum f_h|\,f_h\in\Gamma (U,{\mathbb C}_{W^U_h})\}.
\end{equation}
Let $\mathcal D(W)$ be the sheafification of ${\mathcal D^{pre}(W)}$.  The direct limit of the sheaf $\mathcal D(W)$, as $W$ runs through all reduced effective divisors of $X$, we denote by $\mathcal D(*).$
\begin{remark}\label{ko}  In the remainning of this paper, we will use $\Omega^0$ and $\mathcal O_X$ interchangeably, and use $d\Omega^0$and $\Phi^1$ interchangeably.
\end{remark}

Recall the following lemma due to Hodge and Atiyah.
\begin{lemma}[Lemma $8$ in \cite{HA}]\label{R=D}	
\begin{align*}
&(i)\,R^1(W)\cong\mathcal D(W);\,\,\,\,\\
&(ii)\,R^1(*)\cong\mathcal D(*).\,\,\,\,
\end{align*}
\end{lemma}
\begin{remark}[\cite{HA}]\label{r=d} There is an explicit isomorphism between $R^1(W)$ and $\mathcal D(W)$ as follows. Let $x$ be any point of $X$, and suppose that $W_1\cdots W_t$ are the locally irreducible components of $W$ which pass through $x$, $f_h=0$ being a local equation of $W_h.$ Then $\frac{1}{2\pi \sqrt{-1}}\frac{df_h}{f_h}$ defines an element $r_h$ of $R^1(W)_x.$ The isomorphism is given by
	\begin{displaymath}
	\alpha:\mathcal D(W)\rightarrow R^1(W),\,\,\,
	[1_{W_h}]_x\mapsto r_h.
	\end{displaymath}
	
\end{remark} 

\begin{remark}\label{rrdd} Let $W=\bigcup_{i=1}^l W_i$ be the irreducible decomposition of $W$. Then
$$H^0(X,\mathcal D (W))\cong\bigoplus_{i=1}^l\mathbb C \cdot 1_{W_i},$$
where $1_{W_i}$ is the function taking value 1 in $W_i$ and 0 in $X\backslash W_i$. We also identify the above direct sum with $\bigoplus_{i=1}^l\mathbb C  {W_i}$, the vector space consisting of $\mathbb C$-linear formal sums of divisors $W_1,\cdots,W_l$; call the elements  $\mathbb C$-divisors.
\end{remark}

Hodge and Atiyah considered the following short exact sequences of $\mathbb C$-sheaves:
\begin{equation}\label{Omega}
0\rightarrow\mathbb C\rightarrow\Omega^0(*)\rightarrow d\Omega^0(*)\rightarrow 0\,;
\end{equation}
\begin{equation}\label{R}
0\rightarrow d\Omega^0(*)\rightarrow\Phi^1(*)\rightarrow R^1(*)\rightarrow 0.
\end{equation}
The corresponding long exact sequences of the cohomology groups are
\begin{equation}\label{H1d}
\begin{split}
\rightarrow H^0(X,\Omega^0(*))\rightarrow &H^0(X,d\Omega^0(*))\rightarrow H^1(X,\mathbb C)\rightarrow H^1(X,\Omega^0(*))\rightarrow H^1(X,d\Omega^0(*))\rightarrow \\
&\rightarrow H^2(X,\mathbb C)\rightarrow H^2(X,\Omega^0 (*))\rightarrow\cdots;
\end{split}
\end{equation}
\begin{equation}\label{H0R}
0\rightarrow H^0(X,d\Omega^0(*))\rightarrow H^0(X,\Phi^1(*))\xrightarrow{\rm Res} H^0(X,R^1(*))\xrightarrow{\Delta^0} H^1(X,d\Omega^0(*))\rightarrow\cdots.
\end{equation}
Note that $H^0(X,\Phi^1(*))$ is the vector space of closed meromorphic $1$-forms on $X$.  If under map $\rm Res$ in long exact sequence $(\ref{H0R})$ an element of $H^0(X,\Phi^1(*))$ has  image zero in $H^0(X,R^1(*))$,  we  say that it is a  closed meromorphic $1$-form of the second kind; otherwise we say that it is of the third kind.  It is equivalent to saying that $H^0(X,d\Omega^0(*))$ is the group of  closed meromorphic $1$-forms of the second kind. 
We shall say that a number of $1$-forms are independent if no linear combination of them is equal to the differential of a meromorphic function on $X.$
\medskip

Recall the following well known Serre vanishing theorem.
\begin{theorem}[See Lemma 5, 6 and 7 in \cite{HA}]\label{Kodaira} If $W$ is ample, and if $k$ is a sufficiently large integer, then $H^p(V,\Omega^q(kW))=0$ for $p\geq 1$.
\end{theorem}
It is also proved in \cite{HA} that
\begin{theorem}[Theorem 1 in \S3 of \cite{HA}]\label{cloister}
	The number of independent $1$-forms of the second kind is equal to the first Betti number of $X$.
\end{theorem}

Similarly there are short exact sequences as follows:
\begin{equation}\label{WOmega}
0\rightarrow\mathbb C\rightarrow\Omega^0(*W)\rightarrow d\Omega^0(*W)\rightarrow 0\,;
\end{equation}
\begin{equation}\label{WR}
0\rightarrow d\Omega^0(*W)\rightarrow\Phi^1(*W)\rightarrow R^1(W)\rightarrow 0.
\end{equation}
The corresponding long exact sequences are
\begin{equation}\label{WH1d}
\begin{split}
\rightarrow& H^0(X,\Omega^0(*W))\rightarrow H^0(X,d\Omega^0(*W))\rightarrow H^1(X,\mathbb C)\rightarrow H^1(X,\Omega^0(*W))\rightarrow\\
&\rightarrow H^1(X,d\Omega^0(*W))\xrightarrow{\delta^1} H^2(X,\mathbb C)\rightarrow H^2(X,\Omega^0 (*W))\rightarrow\cdots,
\end{split}
\end{equation}
\begin{equation}\label{WH0R}
0\rightarrow H^0(X,d\Omega^0(*W))\rightarrow H^0(X,\Phi^1(*W))\xrightarrow{\rm Res} H^0(X,R^1(W))\xrightarrow{\Delta^0} H^1(X,d\Omega^0(*W))\rightarrow\cdots.
\end{equation}

In a parallel manner, the following result holds.
\begin{proposition}[Proposition 1 in \S3 of \cite{HA}]\label{*2}
	If $W$ is ample, $H^0(X,d\Omega^0(*W))/\text{Im\,}H^0(X,\Omega^0(*W))$ $\cong H^1(X,\mathbb C)$. Hence a basis for  closed meromorphic $1$-forms of the second kind (modulo differentials of meromorphic functions) can be chosen from forms with singularities on any ample divisor $W$.
\end{proposition}
In fact, we have the following effective version.
\begin{proposition}\label{refined}
	If $H^1(X,\mathcal O_X(kW))=0$ for a certain positive integer $k$, then \begin{equation}
    H^0(X,d\Omega^0(kW))/{\rm Im} H^0(X,\Omega^0(kW))\cong H^1(X,\mathbb C).
	\end{equation}  Moreover, a basis for  closed meromorphic $1$-forms of the second kind (modulo differentials of meromorphic functions) can be chosen from forms with singularities on $W$ and  pole order (at most) $k+1$ along $W$.
\end{proposition}
\noindent{\bf {Proof :}} By definition, we have that
\begin{equation}
0\rightarrow\mathbb C\rightarrow\Omega^0(kW)\rightarrow d\Omega^0(kW)\rightarrow 0.
\end{equation}
The corresponding long exact sequence is:
\begin{equation}
H^0(X,\Omega^0(kW))\rightarrow H^0(X,d\Omega^0(kW))\rightarrow H^1(X,\mathbb C)\rightarrow H^1(X,\Omega^0(kW)).
\end{equation}
Noticing that an element of $H^0(X,d\Omega^0(kW))$ has pole orders no more than $k+1$ along $W$, we conclude Proposition \ref{refined}. $\endpf$
\medskip

We end up this section with the following definitions.

\begin{definition}\label{double} Denote by $\Delta^0$ the homomorphism from $H^0(X,R^1(W))$ to $ H^1(X,d\Omega^0(*W))$  in long exact sequence (\ref{WH0R}); denoted by $\delta^1$ the homomorphism from $H^1(X,d\Omega^0(*W))$ to $H^2(X,\mathbb C)$  in long exact sequence (\ref{WH1d}); call the homomorphism $\delta^1\circ\Delta^0:H^0(X,R^1(W))\rightarrow H^2(X,\mathbb C)$ the double delta map. 
\end{definition}

\begin{definition} For any closed meromorphic $1$-form $\Phi\in H^0(X,\Phi^1(*W))$, we call the image ${\rm Res}$$(\Phi)$ of $\Phi$ under the homomorphism $\rm Res$ the residue divisor of $\Phi.$ Let $W=\bigcup_{i=1}^l W_i$ be the irreducible decomposition of $W$. Noticing Lemma $\ref {R=D}$, Remark $\ref{r=d}$ and Remark $\ref{rrdd}$, $\rm Res$$(\Phi)$ is a $\mathbb C$-linear formal sum of $\{W_i\}_{i=1}^l$ (a $\mathbb C$-divisor$)$ as follows,
	$$Res(\Phi)=\sum_{i=1}^la_iW_i,\,\,{a_i\in\mathbb C}\,\,{\text{for}}\,\,i=1,\cdots,l.$$
\end{definition}
\begin{remark}When $W$ is a normal crossing divisor, we can calculate residue divisor $Res(\Phi)$ by taking terms with the form $\frac{dz_i}{z_i}$ in the Laurent series expansion of $\Phi$ (see the proof of Lemma $\ref{TRU}$).In general, we can calculate $Res(\Phi)$ by taking the contour integrals along small loops around the components of $W$ (see Definition $\ref{tai}$).

\end{remark}

\section{A  geometric interpretation of $\delta^1$ and $\delta^1\circ\Delta^0$ }


Since the original Hodge-Atiyah sequences involve infinite-dimensional cohomology groups, it is not effective to control the pole order and is also abstract for the purpose of a geometric understanding.  In this section, we will introduce a special truncation of short exact sequences $(\ref{WOmega})$ and $(\ref{WR})$. Then by using \v Cech cohomology theory, we are able to derive a geometric interpretation of the homomorphisms $\delta^1$ and $\delta^1\circ\Delta^0$. 
\subsection{Truncation lemma}

\begin{lemma}[Truncation lemma]\label{TRU} Let $X$ be a smooth complex manifold of complex dimension $m$. Let $W$ be a reduced divisor of $X$. There exist short  exact sequences of $\mathbb C$-sheaves  on $X$ as follows:
\begin{equation}\label{iCOME}
0\rightarrow\mathbb C\rightarrow\Omega^0\rightarrow d\Omega^0\rightarrow 0\,;
\end{equation}
\begin{equation}\label{trunc}
0\rightarrow d\Omega^0\xrightarrow{\phi} \Phi^1(W)\xrightarrow{\psi} R^1(W)\rightarrow 0.
\end{equation} 
\end{lemma}
\noindent{\bf {Proof :}} Exact sequence (\ref{iCOME}) is a consequence of Poincar\'e lemma.  To prove exact sequence (\ref{trunc}), it suffices to prove Lemma \ref{TRU} locally at each point $x\in X$. Since  $R^1(W)|_x=0$ and $\phi$ is an isomorphism (by Poincar\'e lemma) when  $x\in X\backslash W$,   Lemma (\ref{trunc}) holds trivially for $x\notin W$.

Now let $x$ be a point of $W$.   Take a neighborhood $U_x$ of $x$ in $X$ such that $W$ is defined by the equation
$f_1f_2\cdots f_l=0,\,\,l\geq 1,$
where $f_1,\cdots,f_l\in\mathcal O(U_x)$ are irreducible holomorphic functions vanishing at $x$ and coprime to each other.  We denote by $W_i$ the zero locus of $f_i$ in $U_x$  for $i=1,\cdots,l.$

Since $(d\Omega^0)_x\subset (\Omega^1)_x$, homomorphism $\phi$ is a well-defined and injective.

Noticing that there is a natural homomorphism from $\Phi^1(W)_x$ to $\Phi^1(*W)_x$ and a homomorphism from $\Phi^1(*W)_x$ to $R^1(W)_x$ (see short exact sequence (\ref{*W})), homomorphism $\psi$ is well-defined. Moreover, $\psi$ is surjective, for 
$R^1(W)_x$ is generated by $\frac{df_1}{f_1},\frac{df_2}{f_2},\cdots,\frac{df_l}{f_l}\in \Phi^1(W)_x$ by Lemma $\ref{R=D}$ and Remark $\ref{r=d}$.

In order to show  complex $(\ref{trunc})$ is exact at the place $\Phi^1(W)_x$, it suffices to prove that if $r\in\Phi^1(W)_x$ and $\psi(r)=0$, then $r$ is the germ of a holomorphic $1$-form.  Notice that if a meromorphic $1$-form has its poles on a subvariety of codimension at least two, then the  $1$-form is actually holomorphic (see Lemma \ref{div}). Therefore, it suffices to prove the exactness at  smooth points of $W$. 

Let $x$ be a smooth point of $W$. Take $r\in\Phi^1(W)_x$ with $\psi(r)=0$. In the following, we denote by $U_t$  for $t>0$ the polydisc $\{(z_1,\cdots,z_m)\big||z_i|<t\,\,{\rm for}\,\,\,i=1,\cdots,m\}$. By a holomoprhic change of coordinates, we can assume that $U_x$ is biholomporphic to $U_1$; $W$ is defined by $z_1=0$ in $U_1$; $r$ takes the form in $U_1$ as
\begin{equation}\label{r}
r=\sum_{p=1}^mr_pdz_p=\sum_{p=1}^m{\Bigg(}\sum_{i_1=-\infty}^{1}\frac{g_{pi_1}(z_{2},\cdots,z_{m})}{z_1^{i_1}}{\Bigg)}dz_p,
\end{equation}
where $g_{pi_1}(z_{2},\cdots,z_m)$ is  holomorphic  in variables $z_{2},\cdots,z_m$ for $p=1,\cdots,m$; for $0<t<1$ and $p=1,\cdots,m$,
\begin{equation}
r_p\cdot z_1=\sum_{i_1=-\infty}^{1}z_1^{1-i_1}{g_{pi_1}(z_{2},\cdots,z_{m})}
\end{equation}is an absolutely convergent series in $U_t$. 
\medskip

{\bf {Claim I :}} $g_{11}(z_{2},\cdots,z_{m})\equiv 0$ in $\widehat U_1$ where $\widehat U_1=\{(z_2,\cdots,z_m)\big||z_i|<1\,\,{\rm for}\,\,\,i=2,\cdots,m\}$.
\medskip

{\bf {Proof of the Claim I :}} Since $\psi(r)=0$ in $U_1$, $r$ locally is a differential of a meromorphic function. Then, each line integral of $r$ along a closed loop in $U_1\backslash\{z_1=0\}$ is zero.  For fixed $(z_2,\cdots,z_m)\in\widehat U_1$, define a loop $\gamma_{z_2\cdots z_m}$ by $(\frac{e^{2\pi \sqrt{-1} t}}{2},z_2,\cdots,z_m)$ for $t\in[0,1]$. Computing the line integral $\int_{\gamma_{z_2\cdots z_m}}r$, we have
\begin{equation}
\begin{split}
\int_{\gamma_{z_2\cdots z_m}}r&=\int_{\gamma_{z_2\cdots z_m}}\Big(\sum_{p=1}^mr_pdz_p\Big)=\int_{\gamma_{z_2\cdots z_m}}\Big(\sum_{p=1}^m\sum_{i_1=-\infty}^{1}\frac{g_{pi_1}(z_{2},\cdots,z_{m})}{z_1^{i_1}}dz_p\Big)\\
&=2\pi \sqrt{-1}\cdot g_{11}(z_2,\cdots,z_m).
\end{split}
\end{equation}
Therefore $g_{11}(z_{2},\cdots,z_{m})\equiv 0$ in $\widehat U_1$. \,\,\,$\endpf$
\medskip

{\bf{Claim II :}}  $g_{p1}(z_{2},\cdots,z_{m})\equiv 0$ in $\widehat U_1$ for $p=2,\cdots,m$.
\medskip

{\bf {Proof of Claim II :}} By Claim I, we can rewrite formula ($\ref{r}$) as 

\begin{equation}
r=\sum_{p=1}^mr_pdz_p=\sum_{p=1}^m{\Bigg(}\sum_{i_1=0}^{\infty}z_1^{i_1}{g_{p(-i_1)}(z_{2},\cdots,z_{m})}{\Bigg)}dz_p+\sum_{p=2}^mz_p^{-1}{g_{p1}(z_{2},\cdots,z_{m})}dz_p.
\end{equation}
Taking the differential of $r$, we have 
\begin{equation}
\begin{split}
0=dr&=\sum_{p=1}^md\Big(\sum_{i_1=0}^{\infty}z_1^{i_1}{g_{p{(-i_1)}}(z_{2},\cdots,z_{m})}\Big)\wedge dz_p+\sum_{p=2}^md{\Bigg(}\frac{g_{p1}(z_{2},\cdots,z_{m})}{z_1^{1}}{\Bigg)}\wedge dz_p\\
&=-\sum_{p=2}^m\frac{g_{p1}(z_{2},\cdots,z_{m})}{z_1^{2}}dz_1\wedge dz_p+\cdots.
\end{split}
\end{equation}
Notice that the coefficient of $dz_1\wedge dz_p$ in $dr$ is 
\begin{equation}
-\frac{g_{p1}(z_{2},\cdots,z_{m})}{z_1^{2}}+h_p(z_1,\cdots,z_m),
\end{equation}
where $h_p$ is a holomorhpic function in $U_1$ and $p=2,\cdots,m$.   Then $g_{p1}(z_{2},\cdots,z_{m})\equiv 0$ in $\widehat U_1$ for $p=2,\cdots,m$. $\endpf$

\medskip
As a conclusion, we proved that $r$ is a holomorphic $1$-form in $U_1$, and hence complex $(\ref{trunc})$ is exact at the place $\Phi^1(W)_x$. We completed the proof of Lemma $\ref{TRU}.$ \,\,\,$\endpf$

\medskip

By a similar argument, we can prove the following truncation with high order poles.
\begin{lemma}[Truncation lemma]\label{kTRU} Let $X$ be a smooth complex manifold of complex dimension $m$. Let $W$ be a reduced divisor of $X$. There exist short  exact sequences of $\mathbb C$-sheaves  on $X$ for $k\geq 1$ as follows:
	\begin{equation}\label{kCOME}
	0\rightarrow\mathbb C\rightarrow\Omega^0(kW)\rightarrow d\Omega^0(kW)\rightarrow 0\,;
	\end{equation}
	\begin{equation}\label{ktrunc}
	0\rightarrow d\Omega^0(kW)\xrightarrow{\phi} \Phi^1((k+1)W)\xrightarrow{\psi} R^1(W)\rightarrow 0.
	\end{equation} 
\end{lemma}
\noindent{\bf Proof of Lemma $\ref{kTRU}$ : } We first prove the following Claim.
\medskip

{\bf Claim :} Let $U$ be a open set of $X$. Let $f\in H^0(U,\mathcal O_X((k+1)W))$ such that $f$ has pole orders at most $k$ in $U$ outside an analytic subvariety of codimension at least two. Then $f\in H^0(U,\mathcal O_X(kW))$.
\medskip

{\bf Proof of Claim :} Noticing that the problem is local, we can assume that $U$ is a complex ball and  $W$ is defined by a reduced holomorphic function $h\in H^0(U,\mathcal O_X)$. Then $f\cdot h^k$ is holomorphic in $U$ outside  an analytic subvariety of codimension at least two; hence, $f\cdot h^k$ is holomorphic in $U$. Since the pole order of $\frac{1}{h}$ on $W$ is $1$, we complete the proof.\,\,\,$\endpf$
\medskip

By the above Claim, we can reduce the problem to a smooth point of $W$. The remaining of the proof is similar to the proof of Lemma \ref{trunc}, and hence we omit it here.\,\,\, $\endpf$

\subsection{Acyclic lemma }
By Lemma $\ref{TRU}$, we have the following short exact sequences of sheaves:
\begin{equation}\label{COME}
0\rightarrow\mathbb C\rightarrow\Omega^0\rightarrow d\Omega^0\rightarrow 0\,.
\end{equation}
We will prove the following acyclic lemma.
\begin{lemma}[Acyclic lemma]\label{AL} Suppose that $X$ is a compact complex manifold. Let $\mathcal U:=\{U_i\}_{i=1}^M$ be a finite, good cover of  $X$ ensured by Lemma $\ref{GC}$ in Appendix I.  Denote by $U_{i_1\cdots i_p}$ the intersection $\bigcap\limits_{j=1}^p U_{i_j}$ for $p\geq 1$ and $1\leq i_1<\cdots<i_p\leq M.$ The following vanishing results hold.
\begin{equation}\label{HC}
H^1(U_{i_1\cdots i_p},\mathbb C)=H^2(U_{i_1\cdots i_p},\mathbb C)=H^3(U_{i_1\cdots i_p},\mathbb C)=\cdots=0,\,p\geq 1;
\end{equation}
\begin{equation}\label{HO}
H^1(U_{i_1\cdots i_p},\Omega^0 )=H^2(U_{i_1\cdots i_p},\Omega^0)=H^3(U_{i_1\cdots i_p},\Omega^0)=\cdots=0,\,p\geq 1;
\end{equation}
\begin{equation}\label{HD}
H^1(U_{i_1\cdots i_p},d\Omega^0)=H^2(U_{i_1\cdots i_p},d\Omega^0)=H^3(U_{i_1\cdots i_p},d\Omega^0)=\cdots=0,\,p\geq 1.
\end{equation}
\end{lemma}

\noindent{\bf{Proof :}} Without loss of generality, we can assume  $U_{i_1\cdots i_p}$ is nonempty; otherwise, the above formulas hold trivially. Since $U_{i_1\cdots i_p}$ is contractible,  the sheaf cohomology group $H^q(U_{i_1\cdots i_p},\mathbb C)$ equals the singular cohomology group $H^q_{sing}(U_{i_1\cdots i_p},\mathbb C)$ for $q\geq 0$ (see \cite{Se} for instance). Since $H^1_{sing}(U_{i_1\cdots i_p},\mathbb C)=H^2_{sing}(U_{i_1\cdots i_p},\mathbb C)=H^3_{sing}(U_{i_1\cdots i_p},\mathbb C)=\cdots=0$, $H^1(U_{i_1\cdots i_p},\mathbb C)=H^2(U_{i_1\cdots i_p},\mathbb C)=H^3(U_{i_1\cdots i_p},\mathbb C)=\cdots=0$. This is formula $(\ref{HC})$.


Since $U_{i_1},\cdots, U_{i_p}$ are Stein, so is $U_{i_1\cdots i_p}$; hence  $H^1(U_{i_1\cdots i_p},\Omega^0)$ $=H^2(U_{i_1\cdots i_p},\Omega^0)=\cdots=0$ by Cartan theorem B, for  $\Omega^0=\mathcal O_X$ is coherent. This is formula $(\ref{HO})$.

Consider the following long exact sequence  associated with short exact sequence $(\ref{COME})$:
\begin{equation*}
\begin{split}
\cdots\rightarrow& H^0(U_{i_1\cdots i_p},\Omega^0)\rightarrow H^0(U_{i_1\cdots i_p},d\Omega^0)\rightarrow H^1(U_{i_1\cdots i_p},\mathbb C)\rightarrow H^1(U_{i_1\cdots i_p},\Omega^0)\rightarrow\\
&\rightarrow H^1(U_{i_1\cdots i_p},d\Omega^0)\rightarrow H^2(U_{i_1\cdots i_p},\mathbb C)\rightarrow H^2(U_{i_1\cdots i_p},\Omega^0)\rightarrow\cdots,
\end{split}
\end{equation*}
Because $H^q(U_{i_1\cdots i_p},\Omega^0)=H^{q+1}(U_{i_1\cdots i_p},\mathbb C)=0$ for $q\geq1$ and $p\geq 1$, we conclude that $H^q(U_{i_1\cdots i_p},d\Omega^0)=0$ for $q\geq1$ and $p\geq 1$. This is formula $(\ref{HD})$. We complete the proof of Lemma \ref{AL}.\,\, $\endpf$
\medskip


Similarly, by Lemma \ref{kTRU} we have  the short exact sequence
	\begin{equation}
	0\rightarrow\mathbb C\rightarrow\Omega^0(kW)\rightarrow d\Omega^0(kW)\rightarrow 0\,,\,\,k\geq 1,
	\end{equation}
and the following lemma.
\begin{lemma}[Acyclic lemma]\label{kAL} Suppose that $X$ is a compact complex manifold. Let $\mathcal U:=\{U_i\}_{i=1}^M$ be a good cover of  $X$ ensured by Lemma $\ref{GC}$ in Appendix I.  The following vanishing results hold for $k\geq 1$:
	\begin{equation}\label{kHO}
	H^1(U_{i_1\cdots i_p},\Omega^0 (kW) )=H^2(U_{i_1\cdots i_p},\Omega^0(kW))=H^3(U_{i_1\cdots i_p},\Omega^0(kW))=\cdots=0,\,p\geq 1;
	\end{equation}
	\begin{equation}\label{kHD}
	H^1(U_{i_1\cdots i_p},d\Omega^0(kW))=H^2(U_{i_1\cdots i_p},d\Omega^0(kW))=H^3(U_{i_1\cdots i_p},d\Omega^0(kW))=\cdots=0,\,p\geq 1.
	\end{equation}
\end{lemma}
\noindent{\bf{Proof :}} The proof is similar to the proof of Lemma \ref{AL} and we omit it here.\,\,$\endpf$

\subsection{The \v Cech cohomology interpretation of homomorphisms}

Let $X$ be a compact complex manifold, $W$ be a reduced divisor of $X$. Then we have the following short exact sequences
\begin{equation}\label{cOME}
0\rightarrow\mathbb C\rightarrow\Omega^0\rightarrow d\Omega^0\rightarrow 0\,;
\end{equation}
\begin{equation}\label{Trunc}
0\rightarrow d\Omega^0\xrightarrow{\phi} \Phi^1(W)\xrightarrow{\psi} R^1(W)\rightarrow 0.
\end{equation} 

Before intepreting the sheaf cohomology as the \v Cech cohomolgoy, we first recall the notions \v Cech resolution and \v Cech complex (see \S 4 of \cite{V} for details).  Let $\mathcal U:=\{U_i\}_{i=1}^M$ be a finite, good cover of $X$. Denote by $j_{*}^{i_1\cdots i_p},p\geq 1,$ the inclusion $U_{i_1\cdots i_p}\xrightarrow{j_{*}^{i_1\cdots i_p}}X$.  Define sheaf $j_{*}^{i_1\cdots i_p}\mathcal F$, for any sheaf of abelian groups $\mathcal F$ on $U_{i_1\cdots i_p}$,  by formula $j_{*}^{i_1\cdots i_p}\mathcal F(V):=\mathcal F(V\cap U_{i_1\cdots i_p})$, where $\mathcal F(\cdot)$ is the notion for the global section functor $\Gamma(\cdot,\mathcal F)$.  Define the sheaf $\mathcal C^k(\mathcal U,\mathcal F)$, for integer $k\geq0$, by $$\mathcal C^k(\mathcal U,\mathcal F):=\bigoplus_{1\leq i_1<\cdots<i_{k+1}\leq M}j_{*}^{i_1\cdots i_{k+1}}\mathcal F.$$ Define the coboundary operator $d:\mathcal F^k\rightarrow\mathcal F^{k+1}$, for integer $k\geq0$, by  formula
$$(d\sigma)_{i_1\cdots i_{k+2}}=\sum\limits_{s}(-1)^{s-1}\sigma_{i_1\cdots\hat i_s\cdots i_{k+2}}|_{V\cap U_{i_1\cdots i_{k+2}}},\,\,1\leq i_1<\cdots <i_{k+2}\leq M,$$
where $\sigma=(\sigma_{j_1\cdots j_{k+1}})$, $\sigma_{j_1\cdots j_{k+1}}\in j_{*}^{j_1\cdots j_{k+1}}\mathcal F(V)=\mathcal F(V\cap U_{j_1\cdots j_{k+1}}),1\leq j_1<\cdots <j_{k+1}\leq M$. One also defines homomorphism $j:\mathcal F\rightarrow C^0(\mathcal U,\mathcal F)$ by $j(\sigma)_i=\sigma|_{V\cap U_i}$ for $\sigma\in\mathcal F(V).$ 
\medskip

We have the following proposition.
\begin{proposition}[Proposition 4.17 in \cite{V}] The (\v Cech) complex
	\begin{equation}\label{resol}
	0\rightarrow\mathcal C^0(\mathcal U,\mathcal F)\xrightarrow{d}\mathcal C^1(\mathcal U,\mathcal F)\xrightarrow{d}\cdots\xrightarrow{d} \mathcal C^n(\mathcal U,\mathcal F)\xrightarrow{d} \mathcal C^{n+1}(\mathcal U,\mathcal F)\xrightarrow{d}\cdots
	\end{equation}
is a resolution of $\mathcal F$. 
\end{proposition}
We call this resolution the \v Cech resolution of $\mathcal F$ associated to the cover $\mathcal U$. 
\begin{definition}Define $\check H ^q(X,\mathcal F)$ to be the $q$th cohomology group of the complex of the global sections
	\begin{equation*}
	C^q(\mathcal U,\mathcal F):=\bigoplus_{1\leq i_1<\cdots<i_{q+1}\leq M}\mathcal F(U_{i_1\cdots i_{q+1}})
	\end{equation*}
of the \v Cech complex $(\ref{resol})$ associated to the cover $\mathcal U$.
\end{definition}
\medskip

In the following, we will consider the \v Cech complexes $(\ref{resol})$ for sheaves $\mathbb C$, $\Omega^0$, $d\Omega^0$, $R^1(W)$ and $\Phi^1(W)$, respectively.
Since the derived functor of the global section functor $\Gamma(X,\cdot)$ is left exact, we have the following isomorphisms:
\begin{equation}\label{check0}
\begin{split}
&H^0(X,\mathbb C)=\check H ^0(X,\mathbb C),\,\,\,\,\,\,\,\,\,\,\,\,\,\,\,\,\,\,\,\,\,\,\,\,\,\,\,\,\,\,\,\,\,\,\,\,\,\,\,\, H^0(X,\Omega^0)=\check H ^0(\mathcal U,\Omega^0),\\
&H^0(X,d\Omega^0)=\check H ^0(\mathcal U,d\Omega^0),\,\,\,\,\,\,\,\,\,\,\,\,\,\,
H^0(X,R^1(W))=\check H ^0(\mathcal U,R^1(W)),\\ &H^0(X,\Phi^1(W))=\check H ^0(\mathcal U,\Phi^1(W)).\\ 
\end{split}
\end{equation}
Moreover, since $H^q(U_{i_1\cdots i_p},\mathbb C)=H^q(U_{i_1\cdots i_p},\Omega^0)$ $=H^q(U_{i_1\cdots i_p},d\Omega^0)$ $=0$ for $p\geq 1$ and $q\geq 1$  by Lemma $\ref{AL}$, there are isomorphisms between the sheaf cohomology groups and the \v Cech cohomology groups (see Theorem 4.41 in \cite{V} for instance) as follows:
\begin{equation}\label{checkp}
\begin{array}{lll}
&H^q(X,\mathbb C)=\check H ^q(\mathcal U,\mathbb C),\,&\,q\geq 1;\\ &H^q(X,\Omega^0)=\check H ^q(\mathcal U,\Omega^0),\,&\,q\geq 1;\\ &H^q(X,d\Omega^0)=\check H ^q(\mathcal U,d\Omega^0),\,&\,q\geq 1.\\ 
\end{array}
\end{equation}

Applying the global section functor and its derived functor to short exact sequence $(\ref{cOME})$ and  long exact sequence $(\ref{resol})$ with $\mathcal F$ $=\mathbb C,$ $\Omega^0$ and $d\Omega^0$,  we have the following commutative diagram, noticing  $H^1(U_{i_1\cdots i_p},\mathbb C)=0$ for $p\geq1$:
\begin{equation}\label{longdo}
\small
\begin{array}{ccccccccc}
0&\rightarrow&\bigoplus\limits_{i_1}\mathbb C(U_{i_1})&\xrightarrow{}&\bigoplus\limits_{i_1} \Omega^0(U_{i_1})&\xrightarrow{\text{}}&\bigoplus\limits_{i_1} d\Omega^0(U_{i_1})&\rightarrow& 0\\
&&\Big\downarrow\rlap{$\scriptstyle d_{\mathbb C}^0$}&&\Big\downarrow\rlap{$\scriptstyle d_{\Omega}^0$}&&\Big\downarrow\rlap{$\scriptstyle d_{d\Omega}^0$}&&\\
0&\rightarrow&\bigoplus\limits_{i_1<i_2}\mathbb C (U_{i_1i_2})&\xrightarrow{}&\bigoplus\limits_{i_1<i_2}\Omega^0 (U_{i_1i_2})&\xrightarrow{D}&\bigoplus\limits_{i_1<i_2} d\Omega^0(U_{i_1i_2})&\rightarrow& 0\\
&&\Big\downarrow\rlap{$\scriptstyle d_{\mathbb C}^1$}&&\Big\downarrow\rlap{$\scriptstyle d_{\Omega}^1$}&&\Big\downarrow\rlap{$\scriptstyle d_{d\Omega}^1$}&&\\
0&\rightarrow&\bigoplus\limits_{i_1<i_2<i_3}\mathbb  C(U_{i_1i_2i_3})&\xrightarrow{}&\bigoplus\limits_{i_1<i_2<i_3}\Omega^0(U_{i_1i_2i_3})&\xrightarrow{\text{}}&\bigoplus\limits_{i_1<i_2<i_3}d\Omega^0(U_{i_1i_2i_3})&\rightarrow& 0\\
&&\Big\downarrow\rlap{$\scriptstyle d_{\mathbb C}^2$}&&\Big\downarrow\rlap{$\scriptstyle d_{\Omega}^2$}&&\Big\downarrow\rlap{$\scriptstyle d_{d\Omega}^2$}&&\\
&&\cdots&&\cdots &&\cdots&&\\
\end{array}.
\end{equation}
Here the horizontal lines of the above commutative diagram are exact. 
By formulas $(\ref{check0})$ and  $(\ref{checkp})$, we conclude that
 $H^0(X,\mathbb C)=\check H^0(\mathcal U,\mathbb C)={\text{Ker}}\, d_{\mathbb C}^0$, $H^1(X,\mathbb C)=\check H^1(\mathcal U,\mathbb C)=\frac{{\text{Ker}}\, d_{\mathbb C}^1}{{\text{Im}}\, d_{\mathbb C}^0}$ and $H^2(X,\mathbb C)=\check H^2(\mathcal U,\mathbb C)=\frac{{\text{Ker}}\, d_{\mathbb C}^2}{{\text{Im}}\, d_{\mathbb C}^1}$;
$H^0(X,\Omega^0)=\check H^0(\mathcal U,\Omega^0)={\text{Ker}}\, d_{\Omega}^0$, $H^1(X,\Omega^0)=\check H^1(\mathcal U,\Omega^0)$ $=\frac{{\text{Ker}}\, d_{\Omega}^1}{{\text{Im}}\, d_{\Omega}^0}$ and $H^2(X,\Omega^0)=\check H^2(\mathcal U,\Omega^0)=\frac{{\text{Ker}}\, d_{\Omega}^2}{{\text{Im}}\, d_{\Omega}^1}$;
$H^0(X,d\Omega^0)=\check H^0(\mathcal U,d\Omega^0)={\text{Ker}}\, d_{d\Omega}^0$ and $H^1(X,d\Omega^0)=\check H^1(\mathcal U,d\Omega^0)=\frac{{\text{Ker}}\, d_{d\Omega}^1}{{\text{Im}}\, d_{d\Omega}^0}$.

Moreover, we also have the following natural commutative diagram:
\begin{equation}\label{do=c}
\small
\begin{array}{ccccccccc}
\rightarrow&\check H^1(\mathcal U,\Omega^0)&\rightarrow & \check H^1(\mathcal U,d\Omega^0)&\xrightarrow{\check\delta^1}&\check H^2(\mathcal U,\mathbb C)&\rightarrow &\check H^2(\mathcal U,\Omega^0)&\rightarrow\\
&\Big\downarrow\rlap{$\scriptstyle \cong$}&&\Big\downarrow\rlap{$\scriptstyle \cong$}&&\Big\downarrow\rlap{$\scriptstyle \cong$}&&\Big\downarrow\rlap{$\scriptstyle \cong$}\\
\rightarrow&H^1(X,\Omega^0)&\rightarrow &H^1(X,d\Omega^0)&\xrightarrow{\delta^1}
&H^2(X,\mathbb C)&\rightarrow&H^2(X,\Omega^0)&\rightarrow\\
\end{array}.
\end{equation}

\bigskip


Noticing $H^1(U_{i_1\cdots i_p},d\Omega^0)=0$ for $p\geq 1$, similarly, we have the following commutative diagram associated with short exact sequence (\ref{Trunc}):
\begin{equation}\label{longr}
\small
\begin{array}{ccccccccc}
0&\rightarrow&\bigoplus\limits_{i_1}d\Omega^0(U_{i_1})&\xrightarrow{}&\bigoplus\limits_{i_1} \Phi^1(W)(U_{i_1})&\xrightarrow{H}&\bigoplus\limits_{i_1} R^1(W)(U_{i_1})&\rightarrow& 0\\
&&\Big\downarrow\rlap{$\scriptstyle d_{d\Omega}^0$}&&\Big\downarrow\rlap{$\scriptstyle d_{\Phi}^0$}&&\Big\downarrow\rlap{$\scriptstyle d_{R}^0$}&&\\
0&\rightarrow&\bigoplus\limits_{i_1<i_2}d\Omega^0(U_{i_1i_2})&\xrightarrow{G}&\bigoplus\limits_{i_1<i_2}\Phi^1(W) (U_{i_1i_2})&\xrightarrow{\text{}}&\bigoplus\limits_{i_1<i_2} R^1(W)(U_{i_1i_2})&\rightarrow& 0\\
&&\Big\downarrow\rlap{$\scriptstyle d_{d\Omega}^1$}&&\Big\downarrow\rlap{$\scriptstyle d_{\Phi}^1$}&&\Big\downarrow\rlap{$\scriptstyle d_{R}^1$}&&\\
0&\rightarrow&\bigoplus\limits_{i_1<i_2<i_3}d\Omega^0(U_{i_1i_2i_3})&\xrightarrow{}&\cdots&\xrightarrow{\text{}}&\cdots&\rightarrow& 0\\
\end{array},
\end{equation}
where the horizontal lines are exact. By formulas $(\ref{check0})$ and  $(\ref{checkp})$, $H^0(X,d\Omega^0)=\check H^0(\mathcal U,d\Omega^0)$ $={\text{Ker}}\, d_{d\Omega}^0$ and $H^1(X,d\Omega^0(kW))=\check H^1(\mathcal U,d\Omega^0)=\frac{{\text{Ker}}\, d_{d\Omega}^1}{{\text{Im}}\, d_{d\Omega}^0}$;
 $H^0(X,\Phi^1(W))=\check H^0(\mathcal U,\Phi^1(W))={\text{Ker}}\, d_{\Phi}^0$; $H^0(X,R^1(W))$ $=\check H^0(\mathcal U,R^1(W))$ $={\text{Ker}}\, d_{R}^0$.

Moreover, we have the following commutative diagram:
\begin{equation}\label{rmapstodo}
\footnotesize
\begin{array}{cccccccc}
\rightarrow \check H^0(\mathcal U,d\Omega^0)&\rightarrow & \check H^0(\mathcal U,\Phi^1(W))&\rightarrow&\check H^0(\mathcal U,R^1(W))&\xrightarrow{\check\Delta^0} &\check H^1(\mathcal U,d\Omega^0)&\rightarrow\cdots\\
\Big\downarrow\rlap{$\scriptstyle \cong$}&&\Big\downarrow\rlap{$\scriptstyle \cong$}&&\Big\downarrow\rlap{$\scriptstyle \cong$}&&\Big\downarrow\rlap{$\scriptstyle \cong$}\\
\rightarrow H^0(X,d\Omega^0)&\rightarrow &H^0(X,\Phi^1(W))&\rightarrow
&H^0(X,R^1(W))&\xrightarrow{\Delta^0}&H^1(X,d\Omega^0)&\rightarrow\cdots\\
\end{array}.
\end{equation}

Combining diagrams $(\ref{do=c})$ and $(\ref{rmapstodo})$, we have a homomorphism $\delta^1\circ\Delta^0$ from $H^0(X,R^1(W))$ to $H^2(X,\mathbb C)$ such that the following diagram is commutative:
\begin{equation}\label{deltaDelta}
\xymatrix{
&H^0(X,R^1(W)) \ar[d]^{\Delta^0} \ar[rd]^{\delta^1\circ\Delta^0}&\\
 0\rightarrow H^0(X,d\Omega^0)\rightarrow H^1(X,\mathbb C)\rightarrow H^1(X,\Omega^0) \ar[r]^{\,\,\,\,\,\,\,\,\,\,\,\,\,\,\,\,\,\,\,\,\,\,\,\,\,\,\,\,\,\,\,\,\,\,\,\,\,\,\,\,\,\,\,\,\,\,\,\,j} &H^1(X,d\Omega^0)\ar[r]^{\,\,\,\,\,\,\delta^1}&H^2(X,\mathbb C).}
\end{equation}

The main result of this section is the following theorem.

\begin{theorem}\label{firstchern}
Let $W$ be a reduced divisor on a compact complex manifold $X$. Let $W=\bigcup_{i=1}^lW_i$ be the irreducible decomposition of $W$. Then the map $\delta^1\circ\Delta^0$ is induced by the first Chern classes as follows,
\begin{equation}\label{chern}
\begin{split}
\delta^1\circ\Delta^0:H^0(X,R^1(W))&\rightarrow H^2(X,\mathbb C),\\
\sum_{i=1}^la_i\cdot1_{W_i}\cong\sum_{i=1}^la_i{W_i}&\mapsto \sum_{i=1}^lc_1(W_i)\otimes_{\mathbb Z}a_i,
\end{split}
\end{equation}
where  $c_1(W_i)$ is the first Chern class of $W_i$ for $i=1,\cdots,l.$  By a slight abuse of notation,  we  call $(\delta^1\circ\Delta^0)(D)$ the first Chern class of $D$ in the  De Rham cohomology for each $D\in H^0(X,R^1(W))$.
\end{theorem}
\begin{remark} Since $H^2(X,\mathbb Z)\otimes_{\mathbb Z}\mathbb C\cong H^2(X,\mathbb C)$, $c_1(W_i)\otimes_{\mathbb Z}1$ can be naturally viewed as an element of $H^2(X,\mathbb C)$. Also recall that, by Lemma $\ref{R=D}$, $$H^0(X,R^1(W))=\bigoplus_{i=1}^lH^0(X,\mathbb C_{W_i})\cong\bigoplus_{i=1}^l\mathbb C\cdot1_{W_i}\cong\bigoplus_{i=1}^l\mathbb C{W_i}.$$.
\end{remark}
\noindent{\bf {Proof :}}  Since $W$ is a reduced divisor, Lemma $\ref{TRU}$ holds; hence short exact sequences (\ref{cOME}) and (\ref{Trunc}) hold. Since $X$ is a compact algebraic manifold, Lemma $\ref{GC}$ and Lemma $\ref{AL}$ hold.  Therefore, we have commutative diagrams (\ref{longdo}),(\ref{do=c}),(\ref{longr}),(\ref{rmapstodo}) and (\ref{deltaDelta}) with respect to a good cover $\mathcal U:=\{U_{i}\}_{i=1}^M$ of $X$.


By the linearity of $\delta^1\circ\Delta^0$, in order to prove formula  (\ref{chern}), it suffices to prove it for $1_{W_1}$.  Notice that, by diagrams $(\ref{do=c})$ and $(\ref{rmapstodo})$, the homomorphism $\delta^1\circ\Delta^0$ is isomorphic to the homomorphism $\check\delta^1\circ\check\Delta^0$ between the corresponding \v Cech cohomology groups.  We will do diagram chasing in diagrams $(\ref{longdo})$ and $(\ref{longr})$ in the following.

Recall that the \v Cech $0$-cocycle of $1_{W_1}$ with respect to $\mathcal U$ is given by $$\bigoplus_{1\leq i_1\leq M}1_{W_1}|_{U_{i_1}}\in\bigoplus\limits_{1\leq i_1\leq M}R^1(W)(U_{i_1})\,\,.$$ 
	
Since $U_i$ is Stein for $i=1,\cdots,M$, $W_1$ is defined by a holomorphic function $f_i=0$ on $U_i$. Define $g_{i_1i_2}:=\frac{f_{i_1}}{f_{i_2}}\in\mathcal O^*(U_{i_1i_2})$ for $i_1,i_2=1,\cdots,M$. Then  $\{g_{i_1i_2}\}$ is a system of transition functions associated with the holomorphic line bundle $W_1$ with respect to $\mathcal U$.

Noticing Remark \ref{r=d} and diagram (\ref{longr}), we define a preimage $\sigma$ of $\bigoplus\limits_{1\leq i_1\leq M}1_{W_1}|_{U_{i_1}}$ under the homomorphism $H:\bigoplus\limits_{1\leq i_1\leq M} \Phi^1(W)(U_{i_1})\xrightarrow{F}\bigoplus\limits_{1\leq i_1\leq M} R^1(W)(U_{i_1})$ by
\begin{equation}
\sigma:=\bigoplus_{1\leq i_1\leq M}\frac{1}{2\pi \sqrt{-1}}\frac{df_{i_1}}{f_{i_1}}\big|_{U_{i_1}}\in\bigoplus\limits_{1\leq i_1\leq M} \Phi^1(W)(U_{i_1}).
\end{equation}
Then $d_{\Phi}^0(\sigma)$ is a \v Cech $1$-cocycle as follows,
\begin{equation}
d_{\Phi}^0(\sigma)=\bigoplus\limits_{1\leq i_1<i_2\leq M}\big(\frac{1}{2\pi \sqrt{-1}}\frac{df_{i_1}}{f_{i_1}}-\frac{1}{2\pi \sqrt{-1}}\frac{df_{i_2}}{f_{i_2}}\big)\big|_{U_{i_1i_2}}\in\bigoplus\limits_{1\leq i_1<i_2\leq M}\Phi^1(W) (U_{i_1i_2}).
\end{equation}


Since $f_{i_1}=g_{i_1i_2}\cdot f_{i_2}$ on $U_{i_1i_2}$ and $g_{i_1i_2}\in\mathcal O^*(U_{i_1i_2})$, we have that
\begin{equation}
\begin{split}
&d_{\Phi}^0(\sigma)(U_{i_1i_2})=\frac{1}{2\pi \sqrt{-1}}\frac{df_{i_1}}{f_{i_1}}-\frac{1}{2\pi \sqrt{-1}}\frac{df_{i_2}}{f_{i_2}}=\frac{1}{2\pi \sqrt{-1}}\frac{d(g_{i_1i_2}f_{i_2})}{(g_{i_1i_2}f_{i_2})}-\frac{1}{2\pi \sqrt{-1}}\frac{df_{i_2}}{f_{i_2}}\\
&=\frac{1}{2\pi \sqrt{-1}}\frac{dg_{i_1i_2}}{g_{i_1i_2}}+\frac{1}{2\pi \sqrt{-1}}\frac{df_{i_2}}{f_{i_2}}-\frac{1}{2\pi \sqrt{-1}}\frac{df_{i_2}}{f_{i_2}}=\frac{1}{2\pi \sqrt{-1}}\frac{dg_{i_1i_2}}{g_{i_1i_2}}=\frac{d(\log g_{i_1i_2})}{2\pi \sqrt{-1}}.
\end{split}
\end{equation}
Define a \v Cech $1$-cocycle $\xi$  by
\begin{equation*}
\xi:=\bigoplus\limits_{1\leq i_1<i_2\leq M}\big(\frac{d(\log g_{i_1i_2})}{2\pi \sqrt{-1}}\big)\big|_{U_{i_1i_2}}\in\bigoplus\limits_{1\leq i_1<i_2\leq M}d\Omega^0(U_{i_1i_2}).
\end{equation*}
It is easy to see that $\xi$ is a preimage of $d_{\Phi}^0(\sigma)$ under the homomorphism $G:\bigoplus\limits_{i_1<i_2}d\Omega^0(U_{i_1i_2})\xrightarrow{G}\bigoplus\limits_{i_1<i_2}\Phi^1(W) (U_{i_1i_2})$.

Fix a base point $a_{i_1i_2}\in U_{i_1i_2}$ for  $1\leq i_1<i_2\leq M.$ Let $\tau$ be a \v Cech $1$-cochain given by
\begin{equation*}
\tau:=\bigoplus\limits_{1\leq i_1<i_2\leq M}\tau_{i_1i_2}|_{U_{i_1i_2}}\in\bigoplus\limits_{1\leq i_1<i_2\leq M}\Omega^0(U_{i_1i_2}),\,\,
\end{equation*}
where
\begin{equation*}
\tau_{i_1i_2}=\int_{a_{i_1i_2}}^z\frac{d(\log g_{i_1i_2})}{2\pi \sqrt{-1}}=\frac{1}{2\pi \sqrt{-1}}\big(\log g_{i_1i_2}(z)-\log g_{i_1i_2}(a_{i_1i_2})\big),
\end{equation*}
and $\log g_{i_1i_2}$ is a branch of the $\log$ function. Then $\tau$ is a preimage of $\xi$ under the homomorphism $D:\bigoplus\limits_{1\leq i_1<i_2\leq M}\Omega^0 (U_{i_1i_2})\xrightarrow{D}\bigoplus\limits_{1\leq i_1<i_2\leq M} d\Omega^0(U_{i_1i_2})$ in diagram $(\ref{longdo})$.

It is clear that $d_{\Omega}^1(\tau)\in\bigoplus\limits_{1\leq i_1<i_2<i_3\leq M}\Omega^0(U_{i_1i_2i_3})$ is a \v Cech $2$-cocycle as follows,
\begin{equation*}
d_{\Omega}^1(\tau)=\bigoplus\limits_{1\leq i_1<i_2<i_3\leq M}d_{\Omega}^1(\tau)_{i_1i_2i_3}\in\bigoplus\limits_{1\leq i_1<i_2<i_3\leq M}\Omega^0(U_{i_1i_2i_3}),\,
\end{equation*}
where
\begin{equation*}
\begin{split}
&d_{\Omega}^1(\tau)_{i_1i_2i_3}=\tau_{i_1i_2}|_{U_{i_1i_2i_3}}-\tau_{i_1i_3}|_{U_{i_1i_2i_3}}+\tau_{i_2i_3}|_{U_{i_1i_2i_3}}=\big\{\frac{1}{2\pi \sqrt{-1}}\big(\log g_{i_1i_2}(z)-\log g_{i_1i_2}(a_{i_1i_2})\big)\\
&-\frac{1}{2\pi \sqrt{-1}}\big(\log g_{i_1i_3}(z)-\log g_{i_1i_3}(a_{i_1i_3})\big)+\frac{1}{2\pi \sqrt{-1}}\big(\log g_{i_2i_3}(z)-\log g_{i_2i_3}(a_{i_2i_3})\big)\big\}\big|_{U_{i_1i_2i_3}}.
\end{split}
\end{equation*}
Since $g_{i_1i_2}\cdot g_{i_2i_3}\cdot g_{i_3i_1}\equiv 1$ on $U_{i_1i_2i_3}$, $d_{\Omega}^1(\tau)_{i_1i_2i_3}$ is a constant function for $z\in U_{i_1i_2i_3}$. Therefore, $d_{\Omega}^1(\tau)\in\bigoplus\limits_{1\leq i_1<i_2<i_3\leq M}\mathbb  C(U_{i_1i_2i_3}).$   Noticing that  $g_{jk}(a_{jk})$ is a constant function defined in $U_{ik}$, we conclude that $d_{\Omega}^1(\tau)$ defines the same  two-cocycle as
\begin{equation}\label{c}
\widetilde{\tau}:=\bigoplus\limits_{1\leq i_1<i_2<i_3\leq M}\frac{1}{2\pi \sqrt{-1}}\big(\log g_{i_1i_2}-\log g_{i_1i_3}+\log g_{i_2i_3}\big)\big|_{U_{i_1i_2i_3}}\in\bigoplus\limits_{1\leq i_1<i_2<i_3\leq M}\mathbb C(U_{i_1i_2i_3}).
\end{equation}
We conclude that $\widetilde\tau$ is the image of $1_{W_1}$ under the map $\check\delta^1\circ\check\Delta^0.$

By Proposition \S 1.1 in \cite{GH}, the above $\widetilde\tau$ corresponds to the first Chern class $c_1(W_1)$ of the divisor $W_1$. We complete the proof of Theorem $\ref{firstchern}$. $\endpf$
\begin{remark}\label{i_2}
	In fact we constructed a homomorphism  $f_1$ from $H^0(X,R^1(W))$ to $H^1(X,\Omega^1)$ such that $\delta^1\circ\Delta^0=i_2\circ f_1$ where $i_2$ is the natrual homomorphism from $H^1(X,\Omega^1)$ to $H^2(X,\mathbb C)$.  By a slight abuse of notaion, we call $f_1(D)$ the first Chern class of $D$ in the Dolbeault cohomology for each $D\in H^0(X,R^1(W))$.
\end{remark}


\subsection{Property (H) and the $\mathcal Q$-flat class of a $\mathbb C$-linear form sum of divisors}
In this subsection, we will introduce the Property $(H)$ of complex manifolds and the $\mathcal Q$-flat classes of holomorphic line bundles and of $\mathbb C$-divisors.
\begin{definition}\label{ph}
	A complex manifold $X$ is said to have Property $(H)$ if $X$ is compact and the following equality holds:
	\begin{equation}
	\dim H^1(X,\mathbb C)=\dim H^{0}(X,d\Omega^0)+\dim H^1(X,\mathcal O_X).
	\end{equation} 
\end{definition}
\begin{remark}\label{nHD}
	Notice that there is a natural long exact sequence of cohomology groups induced by the short exact sequence (\ref{COME}) of $\mathbb C$-sheaves as follows:
	\begin{equation}\label{hd}
	0\rightarrow H^0(X,d\Omega^0)\rightarrow H^1(X,\mathbb C)\rightarrow H^1(X,\mathcal O_X)\rightarrow H^1(X,d\Omega^0)\rightarrow\cdots\,.
	\end{equation}	 
	Hence, in general we only have that $\dim H^1(X,\mathbb C)\leq\dim H^{0}(X,d\Omega^0)+\dim H^1(X,\mathcal O_X)$.
\end{remark}
\begin{remark}\label{kf} Recall that $H^0(X,d\Omega^0)$ is the vector space consisting of the closed holomorphic differential $1$-forms on $X$ and when $X$ is a compact K\"{a}hler manifold (or more generally of Fujiki class $\mathcal C$), each holomorphic $1$-form on $X$ is closed (see \cite{De} or \cite{U}). The following Hodge decomposition for $H^1(X,\mathbb C)$ holds
	\begin{equation}
	0\rightarrow H^0(X,\Omega^1)\rightarrow H^1(X,\mathbb C)\rightarrow H^1(X,\mathcal O_X)\rightarrow 0\,.
	\end{equation}
	Hence, K\"ahler manifolds or manifolds of Fujiki class $\mathcal C$ have Property $(H)$.
\end{remark}
\begin{remark}[See \cite{V}]\label{surf}
	All compact complex surfaces have Property $(H)$.
\end{remark}
\begin{lemma}\label{blow}
	Property $(H)$ is preserved under blow up.
\end{lemma}
\noindent{\bf {Proof of Lemma \ref{blow} :}} Let $f:Y\rightarrow X$ be a blow up map. Since $X$ and $Y$ are smooth, $R^pf_{*}\mathcal O_Y=0$ for $p>0.$ Then $H^1(Y,\mathcal O_Y)=H^1(X,f_{*}\mathcal O_Y)=H^1(X,\mathcal O_X)$. It is easy to verify that $H^0(X,d\Omega^0)=H^0(Y,d\Omega^0)$ and $H^1(X,\mathbb C)=H^1(Y,\mathbb C)$. We complete the proof. \,\,\,$\endpf$
\medskip

Recall that a holomorphic line bundle is said to be flat if its transition functions can be taken as constant functions. We can define an holomorphic invariant for each holomoprhic line bundle, which is called the $\mathcal Q$-flat class, so that a certain integral multiple of the holomorphic line bundle is flat if and only if its $\mathcal Q$-flat class is trivial (see \S 4 of \cite{F})

Since each compact complex manifold has a good cover, we can compute the $\mathcal Q$-flat class map explicitly as follows.
\begin{proposition}[See Theorem 1.8 and Definition 4.3 in \cite{F}] Let $W$ be a holomorphic line bundle over $X$. Let  $\mathcal U=\{U_i\}_{i=1}^M$  be a good cover of $X$ and $\{g_{i_1i_2}\}$ be a system of transition functions  with respect to $\mathcal U$.  Denote by $F$ the homomorphism from $H^1(X,\mathcal O_X^*)$ to $H^1(X,d\Omega^0)$ associating each line bundle with its $\mathcal Q$-flat class. Then $F$ can be computed under the natural isomorphisms $H^1(X,\mathcal O^*)\check H^1(\mathcal U,\mathcal O^*)$ and $H^1(X,d\Omega^0)\cong\check H^1(\mathcal U,d\Omega^0)$ as follows ,
	\begin{equation}
	\begin{split}
	F:\,\,\,\,\,\,\,\,\,\,\,\,H^1(X,\mathcal O^*)&\rightarrow H^1(X,d\Omega^0),\\
	\bigoplus\limits_{1\leq i_1<i_2\leq M}\big(g_{i_1i_2}\big)\big|_{U_{i_1i_2}}
	&\mapsto \bigoplus\limits_{1\leq i_1<i_2\leq M}\big(\frac{d(\log g_{i_1i_2})}{2\pi\sqrt{-1}}\big)\big|_{U_{i_1i_2}}.
	\end{split}
	\end{equation}
\end{proposition}
We call the above homomorphism $F$ the $\mathcal Q$-flat class map.  Recall that there is a natural homormohpism $i_2$ from $H^1(X,\Omega^1)$ to $H^2(X,\mathbb C)$ and a natural homomorphism $j_1$ from $H^1(X,d\Omega^0)$ to $H^1(X,\Omega^1)$. Then we have the following propositions. 
\begin{proposition}[See Theorem 4.5 in \cite{F}] Let $X$ be a compact complex manifold.  The first Chern class maps factor through the $\mathcal Q$-flat class map $F$ as 
	\begin{equation}
	H^1(X,\mathcal O^*)\xrightarrow{F} H^1(X,d\Omega^0)\xrightarrow {j_1} H^1(X,\Omega^1)\xrightarrow{i_2} H^2(X,\mathbb C)\,\,.
	\end{equation}
That is to say, for each holomorphic line bundle $W$ of $X$,  $(j_1\circ F)(W)$ is the first Chern class of $W$ in the Dolbeault cohomology group $H^1(X,\Omega^1)$; $(i_2\circ j_1\circ F)(W)$ is the first Chern class of $W$ in the De Rham cohomology group $H^2(X,\mathbb C)$.
\end{proposition}

We now extend the $\mathcal Q$-flat classes from holomoprhic line bundles to $\mathbb C$-divisors by linearity.
\begin{definition} Suppose $D=\sum_{i=1}^la_i\cdot {W_i}$ where $W_i$ is a divisor of $X$ and $a_i\in\mathbb C$ for $i=1,\cdots,l$.  Define the $\mathcal Q$-flat class of $D$ by $F(D):=\sum_{i=1}^la_i\cdot F(W_i)\in H^1(X,d\Omega^0)$.  
\end{definition}
\begin{remark}
	It is easy to verify that the homomorphism $\Delta^0$ in commutative diagram (\ref{deltaDelta}) is the homomorphism $F$.
\end{remark}
Notice that we have the following lemma.
\begin{lemma}\label{coin}
	 If $X$ has Property $(H)$, then the flat class, the first Chern class in the Dolbeault cohomology and the first Chern class in the De Rham cohomology coincide for each element in $ H^0(X,R^1(W))$.
\end{lemma}
\noindent{\bf {Proof of Lemma \ref{coin} :}} Recalling Remark $\ref{i_2}$ and  commutative diagram $(\ref{deltaDelta})$, we have the following commutative diagram:

\begin{tikzcd}
	&&&&H^0(X,R^1(W)) \ar[d,"\Delta^0"] \ar[dr, "\delta^1\circ\Delta^0"]&\\
	0\ar[r] &H^0(X,d\Omega^0)\ar[r]&H^1(X,\mathbb C)\ar[r]& H^1(X,\Omega^0)\ar[r,"j"]&H^1(X,d\Omega^0)\ar[d,"j_1"]\ar[r,"\delta^1"]&H^2(X,\mathbb C)\\
	&&&&H^1(X,\Omega^1)\ar[ur,"i_2"]&.
\end{tikzcd}
Notice that for each element $D\in H^0(X,R^1(W))$ the $\mathcal Q$-flat class $F(D)=\Delta^0(D)$; the first Chern class in the Dolbeault cohomology is $(j_1\circ\Delta^0)(D)$; the first Chern class in the De Rham cohomology is $(\delta^1\circ\Delta^0)(D)$. Since $X$ has Property $(H)$, the homomorphism $j$ is a zero map, and hence the homomorphisms $\delta^1$, $j_1$ and $i_2$ are injective. Then Lemma \ref{coin} follows.\,\,\,\,$\endpf$
	


\section{Logarithmic forms}
In this section, we will study the logarithmic $1$-forms. First recall the following notions of logarithmic forms (see \cite{No} for instance).

Let $X$ be a compact complex manifold of complex dimension $n$ and  $W$  be  an   effecive  reduced  divisor on $X$.    Fix a point $p\in X$
and   take   irreducible germs of holomorhpic functions $f_j\in\mathcal O_{X,x}$, $1\leq j\leq k$, so that $\{f_1=0\}$, $\cdots$, $\{f_k=0\}$ define   the   local   irreducible   components  of $W$ at $p$. Then  we define  the  sheaf   $\Omega^1(\log W)$ of   germs   of  logarithmic  $1$-forms  along $W$ by
\begin{equation}
\Omega^1_{X,x}(\log W)=\sum_{j=1}^k\mathcal O_{X,x}\frac{df_j}{f_j}+\Omega^1_{X,x}\,\,.
\end{equation}
Notice that the sheaf $\Phi(W)$ is a subsheaf of $\Omega^1(\log W)$.
\subsection{Proofs of Theorems \ref{ri} and \ref{imain}} 

\noindent{\bf {Proof of Theorem \ref{ri} :}} Recall the bottom horizontal line of long exact sequence $(\ref{rmapstodo})$. It is clear the there exists an element $\phi\in H^0(X,\Phi^1(W))$ with residue divisor $D$ if and only if $\Delta^0(D)=0.$ Since $\Delta^0(D)$ is the $\mathcal Q$-flat class of $D$, we conclude Theorem \ref{ri}.\,\,
$\endpf$
\medskip

\noindent{\bf {Proof of Theorem \ref{imain}:}} Recall the commutative diagram $(\ref{deltaDelta})$. When $X$ has Property $(H)$, the homomorphism $j$ is an injection.  Therefore $\Delta^0(D)=0$ if and only if $(\delta^1\circ\Delta^0)(D)=0.$ By Theorem \ref{firstchern} and the Poincar\'e duality, we conclude Theorem \ref{imain}.\,\,
$\endpf$
\medskip

\noindent{\bf {Proof of Corollary \ref{intc}:}} Since the $\mathcal Q$-flat class of a holomorphic line bundle is trivial if only if the line bundle is flat up to some positive multiple (see Theorem 1.11 in \cite{F}), $\Delta_0(W)=\frac{1}{m}\cdot\Delta_0(mW)=0$. By Theorem \ref{ri},  we conclude corollary \ref{intc}.
$\endpf$
\medskip

Noticing Remark \ref{kf} and Remark \ref{surf}, we further have the following corollaries.
\begin{corollary}[The theorem of Weil and Kodaira]
	For a K\"ahler manifold $X$, there exists a closed  logarithmic $1$-form with residue divisor $D$ on $X$ if and only if $D$ is homologous to zero.
\end{corollary}
\begin{corollary}\label{fujiki}
	For a manifold $X$ of Fujiki class $\mathcal C$, there exists a closed logarithmic $1$-form with residue divisor $D$ on $X$ if and only if $D$ is homologous to zero.
\end{corollary}
\begin{remark}
	Corollary \ref{fujiki} can also be derived from Weil and Kodaira's original theorem as follows. Suppose $X$ is a manifold of Fujiki class $\mathcal C$. Then there is a holomorphic map $f$ from a Kahler manifold $X^{\prime}$ to $X$ such that $f$ is generic one to one and surjective. Pull back the residue divisor to $X^{\prime}$ and apply Weil and Kodaira's theorem on $X^{\prime}$. the push forward the form constructed on $X^{\prime}$ to $X$. The resulting form is the desired one.
\end{remark}

\begin{corollary}
	For a compact complex surface $X$, there exists a closed logarithmic $1$-form with residue divisor $D$ on $X$ if and only if $D$ is homologous to zero.
\end{corollary}
\begin{corollary}[Lemma in \S 2.2 of \cite{GH}] Given a finite set of points $\{p_{\lambda}\}$ on compact Riemann surface $S$ and complex numbers $\{a_{\lambda}\}$ such that $\sum a_{\lambda}=0$, there exists a differential of the third kind on $S,$ holomorphic in $S-\{p_{\lambda}\}$ and has residue $a_{\lambda}$ at $p_{\lambda}$.
\end{corollary}

\begin{corollary}\label{general}Let $L$ be a holomorphic line bundle over $X$ and $S_1,S_2,\cdots,S_k$ be generic elements of the complete linear system of $L$. Then $\dim H^0(X,\Phi^1(\widehat S))-\dim H^0(X,\Phi^1)=k-1$ or $k$, where $\widehat S=S_1+\cdots S_k$ is a reduced divisor of $X$.
\end{corollary}
\noindent{\bf {Proof :}} Let $D:=\sum_{i=1}^k{a_i}S_i$ be an element of $H^0(X,R^1(\widehat S))$. Since $S_1,S_2,\cdots,S_k$ are the elements of the complete linear system of $L$, $(\delta^1\circ\Delta^0)(D)=\big(\sum_{i=1}^ka_{i}\big)\cdot c_1(L)$. If $c_1(L)=0$ in $H^2(X,\mathbb C)$, we have that $\dim H^0(X,\Phi^1(\widehat S))-\dim H^0(X,\Phi^1)=k$; otherwise we have that $\dim H^0(X,\Phi^1(\widehat S))-\dim H^0(X,\Phi^1)=k-1$.
$\endpf$
\begin{corollary}\label{ss} Suppose $S$ is a connected, smooth element of an ample line bundle over $X$. Then closed meromorphic $1$-forms with singularities on $S$ are residue free.
\end{corollary}
\noindent{\bf {Proof :}} This is an easy consequence of the fact that $c_1(S)\neq 0.$
$\endpf$

\medskip

In the following, we collect some well known examples on prescribing residues. 

\begin{example}[See \cite{BPV} and \cite{M}]
	Let $X$ be a generic Hopf surface (or generic Hopf manifold). There are finitely many divisors on $X$, each of which is associated with a flat line bundle.  Then we can prescribe  residues on each divisor  (see \cite{P}). On the other hand,  the second singluar cohomology group of $X$ is zero, and hence each divisor is homologous to zero.
\end{example}
\begin{example}[See \cite{N2}]
	Each type VII surfaces has at most finitely many curves. Each curve is homologous to zero and associated with a flat bundle.
\end{example}
\begin{example}[See \cite{N1}]
	Let $X$ be the Iwasawa manifold.  We have that $\dim H^0(X,d\Omega^0)=2$, $\dim H^0(X,\Omega^1)=3$, $\dim H^1(X,\mathcal O_X)=2$ and $\dim H^1(X,\mathbb C)=4$.  Therefore $X$ has property $(H)$ but the Hodge decomposition does not hold (there are non-closed holomorphic $1$-forms on $X$).  On the other hand, $X$ is a fiberation over an abelian variety $T$ of complex dimension 2; the divisors on $X$ are the pull backs of the divisors on $T$.
\end{example}
To end up this section, we would like to ask the following question.
\begin{question}\label{q1}
	Is there  an example of a compact complex manifold $X$ with a $\mathbb C$-divisor $D$ such that $D$ is homologous to zero but there is no closed logarithmic $1$-form  on $X$ with residue $D$?
\end{question}
\subsection{Proof of the decomposition theorem for logarithmic 1-forms}

\noindent{\bf {Proof of Theorem \ref{de} :}} Without loss of generality, we can assume that $W$ is an effective, reduced and normal crossing divisor by Lemma \ref{blow}. In the following, we fix a good cover $\mathcal U=\{U_i\}_{i=1}^M$ of $X$ and choose local coordinates $(z_1^i,\cdots,z^i_n)$ on $U_i$,  for $i=1,\cdots,n$, such that $W\cap U_i=\{z_1\cdots z_{l_i}=0\}$ for certain $0\leq l_i\leq n$ (when $W\cap U_i=\emptyset$, $l_i=0$ by convention).

Let us recall the following construction of the residue divisors for logarithmic $1$-forms (see \cite{B} or \cite{No}). Let $W=\bigcup_{i=1}^mW_i$ be the irreducible decomposition of $W$. Let $\iota_i:\widetilde W_i\rightarrow W_i$ be the normalization of $W_i$ for $i=1,\cdots,m$. Then $\widetilde W:=\coprod_{i=1}^m\widetilde W_i$ is the normalization of $W$ of which the map we denote by $\iota:\widetilde W\rightarrow W$. Noticing that $W$ is a subvariety of $X$, by  abuse of notation, we also denote by $\iota_i$ the map $\widetilde W_i\rightarrow X$, $i=1,\cdots,m$, and $\iota$ the map $\widetilde W\rightarrow X$.  Fix a point $x\in X$. Then  there  is  a   holomorphic  local coordinate system $(z_1,\cdots,z_m)$ in  a  neighborhood $U$ of $x$ such   that $x=(0,\cdots,0)$ and $W\cap U=\{z_1\cdots z_l=0\}\cap U$ where $l$ is an nonegative integer between $0$ and $n$.  Without loss of generality, we can assume that $W_i\cap U=\{z_i=0\}\cap U$ for $i=1,\cdots,l$ and $W_i\cap U=\emptyset$ for $i=l+1,\cdots,m.$ For $\omega\in H^0(X,\Omega^1(\log W))$ we  can write
\begin{equation}
\omega=\sum_{i=1}^l\frac{1}{2\pi \sqrt{-1}}\frac{dz_i}{z_i}\wedge\eta_i+\omega^{\prime}\,\,\,{\rm in}\,\,U,
\end{equation}
where $\eta_i\in H^0(U,\mathcal O_X)$ for $i=1,\cdots,l$ and  $\omega^{\prime}\in H^0(U,\Omega^1)$. Put ${\rm Res}_{\widetilde W_i}(\omega)=\iota^*_i(\eta)$ in $\widetilde W_i\bigcap \iota^{-1}(U)$ for $i=1,\cdots,m$. Then ${\rm Res}_{\widetilde W_i}(\omega)$ is  	globally well-defined and 
\begin{equation}
{\rm Res}_{\widetilde W_i}(\omega)\in H^0(\widetilde W_i,\mathcal O_{\widetilde W_i})\,\,\,{\rm for\,\,}i=1,\cdots,m.
\end{equation}
By pushing forward, we have the following short exact sequence of $\mathcal O_X$-sheaves on $X$
\begin{equation}\label{loga}
0\rightarrow \Omega^1_X\rightarrow \Omega^1_X(\log W)\xrightarrow{{\rm Res}}\iota_*\mathcal O_{\widetilde W}\rightarrow 0\,\,.
\end{equation}

We next consider the \v Cech cohomology groups associated with exact sequence $(\ref{loga})$ with respect to $\mathcal U$.  In the same way as in \S 3.3, the following commutative diagram holds
\begin{equation}\label{logar}
\small
\begin{array}{ccccccccc}
0&\rightarrow&\bigoplus\limits_{i_1}\Omega_X^1(U_{i_1})&\xrightarrow{}&\bigoplus\limits_{i_1} \Omega_X^1(\log W)(U_{i_1})&\xrightarrow{H}&\bigoplus\limits_{i_1} \iota_*\mathcal O_{\widetilde W}(U_{i_1})&\rightarrow& 0\\
&&\Big\downarrow\rlap{$\scriptstyle d_{\Omega}^0$}&&\Big\downarrow\rlap{$\scriptstyle d_{\log}^0$}&&\Big\downarrow\rlap{$\scriptstyle d_{W}^0$}&&\\
0&\rightarrow&\bigoplus\limits_{i_1<i_2}\Omega_X^1(U_{i_1i_2})&\xrightarrow{G}&\bigoplus\limits_{i_1<i_2}\Omega_X^1(\log W) (U_{i_1i_2})&\xrightarrow{\text{}}&\bigoplus\limits_{i_1<i_2} \iota_*\mathcal O_{\widetilde W}(U_{i_1i_2})&\rightarrow& 0\\
&&\Big\downarrow\rlap{$\scriptstyle d_{\Omega}^1$}&&\Big\downarrow\rlap{$\scriptstyle d_{\log}^1$}&&\Big\downarrow\rlap{$\scriptstyle d_{W}^1$}&&\\
0&\rightarrow&\bigoplus\limits_{i_1<i_2<i_3}\Omega_X^1(U_{i_1i_2i_3})&\xrightarrow{}&\cdots&\xrightarrow{\text{}}&\cdots&\rightarrow& 0\\
\end{array}.
\end{equation}
Moreover, we have
\begin{equation}\label{logarx}
\footnotesize
\begin{array}{cccccccc}
\check H^0(\mathcal U,\Omega_X^1)&\rightarrow & \check H^0(\mathcal U,\Omega_X^1(\log W))&\xrightarrow{\rm Res}&\check H^0(\mathcal U,\iota_*\mathcal O_{\widetilde W})&\xrightarrow{\check\Delta^0} &\check H^1(\mathcal U,\Omega_X^1)&\rightarrow\cdots\\
\Big\downarrow\rlap{$\scriptstyle \cong$}&&\Big\downarrow\rlap{$\scriptstyle \cong$}&&\Big\downarrow\rlap{$\scriptstyle \cong$}&&\Big\downarrow\rlap{$\scriptstyle \cong$}\\
H^0(X,\Omega_X^1)&\rightarrow &H^0(X,\Omega_X^1(\log W))&\xrightarrow{\rm Res}
&H^0(X,\iota_*\mathcal O_{\widetilde W})&\xrightarrow{\Delta^0}&H^1(X,\Omega_X^1)&\rightarrow\cdots\\
\end{array}.
\end{equation}
Notice that
\begin{equation}
H^0(X,\iota_*\mathcal O_{\widetilde W})\cong H^0(\widetilde W,\mathcal O_{\widetilde W})\cong H^0(W,\mathcal O_{W})\cong\bigoplus_{i=1}^m\mathbb C \cdot 1_{W_i}\cong H^0(X,R^1(W))\,\,.
\end{equation}
Hence each element $\sigma\in\check H^0(\mathcal U,\iota_*\mathcal O_{\widetilde W})$ can be represented by a \v Cech 0-cocycle
\begin{equation}
\sigma=\bigoplus_{1\leq i_1\leq M}(\sum_{k=1}^ma_k\cdot 1_{W_k})\big|_{U_{i_1}}\in\bigoplus\limits_{1\leq i_1\leq M} \iota_*\mathcal O_{\widetilde W}(U_{i_1})\,\,,
\end{equation}
where $a_k\in\mathbb C$ for $k=1,\cdots,m$.

We will show that the homomorphism $\Delta^0$ between $ H^0(X,\iota_*\mathcal O_{\widetilde W})$ and $ H^1(X,\Omega_X^1)$ is induced by the first Chern classes as follows.
\medskip

\noindent{\bf Claim :} The map $\Delta^0:H^0(X,\iota_*\mathcal O_{\widetilde W})\rightarrow H^1(X,\Omega_X^1)$ is given by
\begin{equation}
\Delta^0:\sum_{j=1}^ma_j\cdot W_j\mapsto \sum_{j=1}^ma_j\cdot c_1(W_j),
\end{equation}
where $c_1(W_j)$ is the first Chern class of the holomorphic line bundle $[W_j]$ in the Dolbeault cohomology group $H^1(X,\Omega^1_X)$. By abuse of notation, we call $\Delta^0(D)$ the first Chern class of $D$ in the Dolbeault cohomology for each $D\in H^0(X,\iota_*\mathcal O_{\widetilde W})$.
\medskip

\noindent {\bf Proof of Claim :} Since $\Delta^0$ is a linear map, without loss of generality, it suffices to prove Claim when $a_1=1,a_2=a_3=\cdots=a_m=0$. Suppose $W_1$ is defined by a holomorphic function $f_i$ in $U_i$ for $i=1,\cdots,M$ and the transition function $g_{ij}=\frac{f_i}{f_j}$ in $U_{ij}$ for $i,j=1,\cdots,M$.
Then the element $1\cdot W_1\in H^0(X,\iota_*\mathcal O_{\widetilde W}) $ can be represented by a \v Cech 0-cocycle $\sigma$ as
\begin{equation}
\sigma=\bigoplus_{1\leq i_1\leq M}( 1_{W_1})\big|_{U_{i_1}}\in\bigoplus\limits_{1\leq i_1\leq M} \iota_*\mathcal O_{\widetilde W}(U_{i_1})\,\,.
\end{equation}
A preimage $\eta$ of $\sigma$ under $H$ can be taken as
\begin{equation}
\eta=\bigoplus_{1\leq i_1\leq M}\big(\frac{1}{2\pi \sqrt{-1}}\frac{df_{i_1}}{f_{i_1}}\big)\big|_{U_{i_1}}\in\bigoplus\limits_{1\leq i_1\leq M}\Omega_X^1(\log W)(U_{i_1})\,\,.
\end{equation}
The \v Cech one-cocycle $d_{\log}^0(\eta)$ takes the form
\begin{equation}
d_{\log}^0(\eta)=\bigoplus\limits_{1\leq i_1<i_2\leq M}\big(\frac{1}{2\pi \sqrt{-1}}\frac{df_{i_1}}{f_{i_1}}-\frac{1}{2\pi \sqrt{-1}}\frac{df_{i_2}}{f_{i_2}}\big)\big|_{U_{i_1i_2}}\in\bigoplus\limits_{1\leq i_1<i_2\leq M}\Omega_X^1(\log W) (U_{i_1i_2}).
\end{equation}
Since $f_{i_1}=g_{i_1i_2}f_{i_2}$, 
\begin{equation}
\frac{1}{2\pi \sqrt{-1}}\frac{df_{i_1}}{f_{i_1}}-\frac{1}{2\pi \sqrt{-1}}\frac{df_{i_2}}{f_{i_2}}=\frac{1}{2\pi \sqrt{-1}}{d(\log g_{i_1i_2})}\in\Omega_X^1(U_{i_1i_2}).
\end{equation}
Therefore we can lift $\eta$ to a \v Cech 1-cocycle $\xi$ as
\begin{equation}
\xi:=\bigoplus\limits_{1\leq i_1<i_2\leq M}\big(\frac{1}{2\pi \sqrt{-1}}{d(\log g_{i_1i_2})}\big)\big|_{U_{i_1i_2}}\in\bigoplus\limits_{1\leq i_1<i_2\leq M}\Omega_X^1 (U_{i_1i_2}).
\end{equation}
By Proposition \S 1.1 in \cite{GH}, we conclude that $\xi$ is the first Chern class of $W_1$ as an $(1,1)$ form. Therefore, we complete the proof of Claim.\,\,$\endpf$ 
\medskip

Now we proceed to prove Theorem \ref{de}. Let $\omega\in H^0(X,\Omega_X^1(\log W))$. Then the first Chern class of Res$(\omega)$ in the Dolbeault cohomology is trivial. Siince $X$ has Property $(H)$, the $\mathcal Q$-flat class of Res$(\omega)$ is trivial by Remark \ref{coin}. Therefore, there is a closed meromorphic $1$-form $\omega_1$ with the residue class Res$(\omega)$. Since the residue divisor of $\omega-\omega_1$  is zero, $\omega-\omega_1$ is a holomorphic $1$-form by long exact sequence $(\ref{logarx})$. We finish the proof of Theorem \ref{de}.\,\,\,\,\,\,\,$\endpf$
\begin{remark}
	Theorem \ref{de} holds under a slightly weaker condition that in the  Fr\"olicher spectral sequence of $X$,
	\begin{equation}
	E_2^{0,1}=E_{3}^{0,1}\,\,.
	\end{equation}
\end{remark}

\begin{example}[See \cite{N1}]
	Let $X$ be an Iwasawa manifold.  Then $X$ has property $(H)$ but does not have the Hodge decomposition. There are non-closed holomorphic $1$-forms on $X$.
\end{example}


\section{Forms with high order poles}
In this section, we will study the closed meromorphic $1$-froms with poles of higher order.  We first refine Hodge and Atiyah's criterion in the following way.
\begin{theorem}\label{knc}
Let X be a compact complex manifold and $W$ be a reduced divisor of $X$.   Let $D\in H^0(X,R^1(W))$ be a $\mathbb C$-linear formal sum of divisors of $X$. Let $k$ be an nonnegative integer. Then the following statements are equivalent:
\begin{enumerate}[]
	\item $\bullet$  $J_k(F(D))$ is trivial in $H^1(X,d\Omega^0(kW))$, where $J_k:H^1(X,d\Omega^0)\rightarrow H^1(X,d\Omega^0(kW))$ is the natural homomorphism associated with the homomorphism of sheaves $j_k:d\Omega^0\rightarrow d\Omega^0(kW)$;
	\item $\bullet$  there is a closed logarithmic $1$-form $\phi\in H^0(X,\Phi((k+1)W))$ with residue divisor $D$. 
\end{enumerate}
\end{theorem}
\noindent{\bf Proof of Theorem \ref{knc} :}
By Lemma \ref{kTRU}, we have the following short exact sequence of sheaves:
\begin{equation}\label{kktrunc}
0\rightarrow d\Omega^0(kW)\xrightarrow{\phi} \Phi^1((k+1)W)\xrightarrow{\psi} R^1(W)\rightarrow 0.
\end{equation} 
By Lemma \ref{kAL} we can compute the sheaf cohomology groups by \v Cech cohomology groups in the same as the approach in  \S 3. In particular,  we have the following commutative diagram:

\begin{equation}\label{klongr}
\small
\begin{array}{ccccccccc}
0&\rightarrow&\bigoplus\limits_{i_1}d\Omega^0(kW)(U_{i_1})&\xrightarrow{}&\bigoplus\limits_{i_1} \Phi^1((k+1)W)(U_{i_1})&\xrightarrow{}&\bigoplus\limits_{i_1} R^1(W)(U_{i_1})&\rightarrow& 0\\
&&\Big\downarrow\rlap{$\scriptstyle d_{d\Omega}^0$}&&\Big\downarrow\rlap{$\scriptstyle d_{\Phi}^0$}&&\Big\downarrow\rlap{$\scriptstyle d_{R}^0$}&&\\
0&\rightarrow&\bigoplus\limits_{i_1<i_2}d\Omega^0(kW)(U_{i_1i_2})&\xrightarrow{}&\bigoplus\limits_{i_1<i_2}\Phi^1((k+1)W) (U_{i_1i_2})&\xrightarrow{\text{}}&\bigoplus\limits_{i_1<i_2} R^1(W)(U_{i_1i_2})&\rightarrow& 0\\
&&\Big\downarrow\rlap{$\scriptstyle d_{d\Omega}^1$}&&\Big\downarrow\rlap{$\scriptstyle d_{\Phi}^1$}&&\Big\downarrow\rlap{$\scriptstyle d_{R}^1$}&&\\
0&\rightarrow&\bigoplus\limits_{i_1<i_2<i_3}d\Omega^0(kW)(U_{i_1i_2i_3})&\xrightarrow{}&\cdots&\xrightarrow{\text{}}&\cdots&\rightarrow& 0\\
\end{array}.
\end{equation}

Moreover, we also have the following commutative diagrams:
\begin{equation}\label{krmapstodo}
\footnotesize
\begin{array}{cccccccc}
\check H^0(\mathcal U,d\Omega^0(kW))&\rightarrow & \check H^0(\mathcal U,\Phi^1((k+1)W))&\rightarrow&\check H^0(\mathcal U,R^1(W))&\xrightarrow{\check\Delta^0_k} &\check H^1(\mathcal U,d\Omega^0(kW))&\rightarrow\cdots\\
\Big\downarrow\rlap{$\scriptstyle \cong$}&&\Big\downarrow\rlap{$\scriptstyle \cong$}&&\Big\downarrow\rlap{$\scriptstyle \cong$}&&\Big\downarrow\rlap{$\scriptstyle \cong$}\\
H^0(X,d\Omega^0(kW))&\rightarrow &H^0(X,\Phi^1((k+1)W))&\rightarrow
&H^0(X,R^1(W))&\xrightarrow{\Delta^0_k}&H^1(X,d\Omega^0(kW))&\rightarrow\cdots\\
\end{array}.
\end{equation}

For homomorphism $j_k:d\Omega^0\rightarrow d\Omega^0(kW)$, we have a natural homorphism between the \v Cech complexes as 
\begin{equation}\label{kd}
\small
\begin{array}{ccc}
\bigoplus\limits_{i_1}d\Omega^0(U_{i_1})&\xrightarrow{}&\bigoplus\limits_{i_1} d\Omega^0(kW)(U_{i_1})\\
\Big\downarrow\rlap{$\scriptstyle d_{d\Omega}^0$}&&\Big\downarrow\rlap{$\scriptstyle d_{\Phi}^0$}\\
\bigoplus\limits_{i_1<i_2}d\Omega^0(U_{i_1i_2})&\xrightarrow{}&\bigoplus\limits_{i_1<i_2}d\Omega^0(kW)(U_{i_1i_2})\\
\Big\downarrow\rlap{$\scriptstyle d_{d\Omega}^1$}&&\Big\downarrow\rlap{$\scriptstyle d_{\Phi}^1$}\\
\bigoplus\limits_{i_1<i_2<i_3}d\Omega^0(U_{i_1i_2i_3})&\xrightarrow{}&\bigoplus\limits_{i_1<i_2<i_3}d\Omega^0(kW)(U_{i_1i_2i_3})\\
\end{array}.
\end{equation}
Hence, the following commutative diagram holds:
\begin{equation}\label{kj}
\footnotesize
\begin{array}{cccccccc}
\check H^1(\mathcal U,d\Omega^0)&\xrightarrow{\check J_k} & \check H^1(\mathcal U,d\Omega^0(kW))\\
\Big\downarrow\rlap{$\scriptstyle \cong$}&&\Big\downarrow\rlap{$\scriptstyle \cong$}\\
H^1(X,d\Omega^0)&\xrightarrow{J_k} &H^1(X,d\Omega^0(kW))\\
\end{array},
\end{equation}
where $\check J_k$ is induced from diagram (\ref{kd}).
Repeating the diagram chasing, we have that the map $\Delta_k^0$ factors through the $\mathcal Q$-flat class map $F$, that is, 
\begin{equation}\label{fact}
H^0(X,R^1(W))\xrightarrow{F}H^1(X,d\Omega^0)\xrightarrow{J_k}H^1(X,d\Omega^0(kW)
\end{equation}
where $J_k\circ F=\Delta_k^0$.

Therefore, by the exactness of diagram (\ref{krmapstodo}) we complete the proof of Theorem \ref{knc}. \,\,\,$\endpf$
\medskip

\noindent{\bf Proof of Theorem \ref{nc} :} First notice that for each closed meromorphic $1$-form $\phi\in H^0(X,\Phi(*))$, there is a reduced divisor $W$ of $X$ and an nonnegative integer $k$ such that $\phi\in H^0(X,\Phi(kW))$. Consider the  short exact sequence
\begin{equation}\label{kkCOME}
0\rightarrow\mathbb C\rightarrow\Omega^0(kW)\rightarrow d\Omega^0(kW)\rightarrow 0\,.
\end{equation}
	
In addition to commutative diagram (\ref{klongr}), we have the following commutative diagrams: 
\begin{equation}\label{klongdo}
\small
\begin{array}{ccccccccc}
0&\rightarrow&\bigoplus\limits_{i_1}\mathbb C(U_{i_1})&\xrightarrow{}&\bigoplus\limits_{i_1} \Omega^0(kW)(U_{i_1})&\xrightarrow{\text{}}&\bigoplus\limits_{i_1} d\Omega^0(kW)(U_{i_1})&\rightarrow& 0\\
&&\Big\downarrow\rlap{$\scriptstyle d_{\mathbb C}^0$}&&\Big\downarrow\rlap{$\scriptstyle d_{\Omega}^0$}&&\Big\downarrow\rlap{$\scriptstyle d_{d\Omega}^0$}&&\\
0&\rightarrow&\bigoplus\limits_{i_1<i_2}\mathbb C (U_{i_1i_2})&\xrightarrow{}&\bigoplus\limits_{i_1<i_2}\Omega^0(kW) (U_{i_1i_2})&\xrightarrow{D}&\bigoplus\limits_{i_1<i_2} d\Omega^0(kW)(U_{i_1i_2})&\rightarrow& 0\\
&&\Big\downarrow\rlap{$\scriptstyle d_{\mathbb C}^1$}&&\Big\downarrow\rlap{$\scriptstyle d_{\Omega}^1$}&&\Big\downarrow\rlap{$\scriptstyle d_{d\Omega}^1$}&&\\
0&\rightarrow&\bigoplus\limits_{i_1<i_2<i_3}\mathbb  C(U_{i_1i_2i_3})&\xrightarrow{}&\bigoplus\limits_{i_1<i_2<i_3}\Omega^0(kW)(U_{i_1i_2i_3})&\xrightarrow{\text{}}&\bigoplus\limits_{i_1<i_2<i_3}d\Omega^0(kW)(U_{i_1i_2i_3})&\rightarrow& 0\\
&&\Big\downarrow\rlap{$\scriptstyle d_{\mathbb C}^2$}&&\Big\downarrow\rlap{$\scriptstyle d_{\Omega}^2$}&&\Big\downarrow\rlap{$\scriptstyle d_{d\Omega}^2$}&&\\
&&\cdots&&\cdots &&\cdots&&\\
\end{array}
\end{equation}
and
\begin{equation}\label{kdo=c}
\small
\begin{array}{ccccccccc}
\rightarrow&\check H^1(\mathcal U,\Omega^0(kW))&\rightarrow & \check H^1(\mathcal U,d\Omega^0(kW))&\xrightarrow{\check\delta^1_k}&\check H^2(\mathcal U,\mathbb C)&\rightarrow &\check H^2(\mathcal U,\Omega^0(kW))&\rightarrow\\
&\Big\downarrow\rlap{$\scriptstyle \cong$}&&\Big\downarrow\rlap{$\scriptstyle \cong$}&&\Big\downarrow\rlap{$\scriptstyle \cong$}&&\Big\downarrow\rlap{$\scriptstyle \cong$}\\
\rightarrow&H^1(X,\Omega^0(kW))&\rightarrow &H^1(X,d\Omega^0(kW))&\xrightarrow{\delta^1_k}
&H^2(X,\mathbb C)&\rightarrow&H^2(X,\Omega^0(kW))&\rightarrow\\
\end{array}.
\end{equation}

Combining diagrams $(\ref{krmapstodo})$ and $(\ref{kdo=c})$, we have a homomorphism $\delta^1_k\circ\Delta^0_k$ from $H^0(X,R^1(W))$ to $H^2(X,\mathbb C)$ such that the following commutative diagram holds:
\begin{equation}\label{kdeltaDelta}
\xymatrix{
	&H^0(X,R^1(W)) \ar[d]^{\Delta^0_k} \ar[rd]^{\delta^1_k\circ\Delta^0_k}&&\\
	H^1(X,\Omega^0(kW))\ar[r] &H^1(X,d\Omega^0(kW))\ar[r]^{\,\,\,\,\,\,\,\,\,\,\,\,\delta^1_k}&H^2(X,\mathbb C)\ar[r] &H^2(X,\Omega^0(kW)).}
\end{equation}
Moreover, by diagram chasing, we can conclude that $(\delta_k^1\circ\Delta^0_k)(D)$ is the first Chern class of $D$ in the De Rham cohomology for each $\mathbb C$-divisor $D$.  By Poincar\'e duality, we complete the proof of Theorem \ref{nc}.\,\,\,\,\,$\endpf$
\medskip

Next we show that for manifolds with Property ($H$) a further decomposition is possible.
\medskip

\noindent{\bf Proof of Theorem \ref{df} :} Let $\phi\in H^0(X,\Phi(*))$. Choose a reduced divisor $W$ and a nonnegative integer $k$ such that $\phi\in H^0(X,\Phi(kW))$. Recalling commutative diagram (\ref{kdeltaDelta}) and the fact that $\Delta^0_k$ factors through $H^1(X,d\Omega^0)$, we have the following commutative diagram:
\begin{equation}\label{tower}
\xymatrix{
	H^0(X,R^1(W)) \ar[d]^{F} \ar[rdd]^{c_1}&\\
	H^1(X,d\Omega) \ar[d]^{J_k} \ar[rd]^{\delta^1}&\\
	H^1(X,d\Omega^0(kW))\ar[r]^{\,\,\,\,\,\,\,\,\,\,\,\,\delta^1_k}&H^2(X,\mathbb C). }
\end{equation}
Since $X$ has Property ($H$), $\delta^1$ is injective; hence $F({\rm Res}(\phi))=0$. Therefore, we can find a closed logarithmic form $\widetilde\phi\in H^0(X,\Phi(W))$ such that ${\rm Res}(\phi)={\rm Res}(\widetilde\phi)$. Then $\psi:=\phi-\widetilde\phi$ is of the second kind. We complete the proof of Corollary \ref{df}\,\,\,\,$\endpf$

\section{Constructing pluriharmonic functions with mild poles}
In this section, we will investigate Question \ref{Q2}. Our construction of pluriharmonic functions is by integrating closed meromorphic $1$-forms. Notice that the integration often results in a multivalued function; the obstructions are twofold, that is, the long period vectors and the short period vectors (to de defined in the following).  In order to get a single-valued function,  we construct  a conjugate, namely, a closed anti-meromorphic $1$-form, for each closed meromorphic $1$-form, so that the cancellation of the periods is possible when integrating the sum of the pair. 

In this following, we always assume that $X$ is a compact algebraic manifold.
 
\subsection{Definition and properties of gardens, period vectors and  pairs}
\begin{definition}\label{garden} 
Let $X$ be a compact algebraic manifold. 	Let $W$ be a reduced divisor of $X$ with the irreducible decomposition $W=\bigcup_{j=1}^lW_j$.  Fix a point $p\in X\backslash W.$ Let $\{\gamma_1,\cdots,\gamma_m\}$ be a basis of $H_1(X,\mathbb C)$ such that $\gamma_i$
is a smooth Jordan curve based at $p$ and contained in $X\backslash W$ for $i=1,\cdots,m$.	We call the quadruple $\big(X,W,(W_1,\cdots,W_l),(\gamma_1,\cdots,\gamma_m)\big)$ a garden.
\end{definition}


\begin{definition}\label{long} Suppose $A=\big(X,W,(W_1,\cdots,W_l),(\gamma_1,\cdots,\gamma_m)\big)$ is a garden. 
Let $\Phi$ be a closed meromorphic $1$-form  on $X$ with singularities on $W$.   We call $(b_1,\cdots,b_m)$ the long period vector of $\Phi$ with respect to $A$ if
$$b_i=\int_{\gamma_i}\Phi\,\,\,{\text{for}}\,\,i=1,\cdots,m.$$
\end{definition}

Then we have the following lemma.
\begin{lemma}\label{connect} Let $W=\bigcup_{i=1}^lW_{i}$ be the irreducible decomposition of reduced divisor $W.$ Let $\Phi$ be a closed meromorphic $1$-form on $X$ with singularities on $W$. Let  $a_1,a_2$ be two smooth points of $W_i$; suppose  $D_1,D_2\subset X$ are two one-dimensional analytic discs intersecting $W_i$ transversally at $a_1$, $a_2$, respectively; suppose $S_k$ is a small disc contained in $D_k$ centered at $a_k$ with $\partial S_k$ the counterclockwise oriented circle boundary of $S_k$ for $k=1,2.$  Then the following integrals are equal,
	\begin{equation*}
	\int_{\partial S_1}\Phi=\int_{\partial S_2}\Phi.
	\end{equation*}
\end{lemma}
\noindent{\bf {Proof:}} Since $W_i$ is irreducible, there is a proper subvariety ${\text{Sing}\,W_i}$ of $W_i$, such that $a_1,a_2\in W_i^{\text{Reg}}:=W_i\backslash {\text{Sing}\,W_i}$ and $W_i^{\text{Reg}}$ is connected. Hence, we can find a smooth curve $\gamma\subset W_i^{\text{Reg}}$ connecting $a_1,a_2$, and an open set $U$ of $\gamma$ in $W_i^{\text{Reg}}$ such that $\gamma\subset U\subset\subset W_i^{\text{Reg}}$.  By Theorem (5.2) in \S 4 of \cite{Hi},  $U$ has a tubular neighborhood $E$  in $X$ in the sense that 
\begin{equation*}
\begin{split}
\pi:E\rightarrow U \,\,{\text{is a disc bundle}};\,\,
J:E\hookrightarrow X\,\,{\text{is an embedding}}.
\end{split}
\end{equation*}
Notice that $J(E)$ is an open set of $X$ and the zero section of $E$ is mapped diffeomorphically to $U.$ Without loss of generality, we assume $J(E)\cap W=U.$ 

When $S_i$ is small enough, we have that $S_i\subset J(E)$ and $\partial S_i$ is homotopic to the counterclockwise oriented boundary $\Gamma_i$ of $J(\pi^{-1}(a_i))$ for $i=1,2$. Notice that  the homotopies can be choosen away from $W$;  moreover, $\Gamma_1$ and $\Gamma_2$ are diffeomorphic in the circle bundle (the boundary of $E$) over $\gamma$.  Therefore, in a same way as the proof of Lemma \ref{smoothing} in Appendix II, we can find smooth two-chains $\sum_{i=1}^pa_i\sigma_i^1$, $\sum_{j=1}^qb_j\sigma_j^2$ and $\sum_{k=1}^rc_k\sigma_k^3$, such that they are disjoint from $W$ and
$$\partial(\sum_{i=1}^pa_i\sigma_i^1)=\partial S_1-\Gamma_1;\,\,\partial(\sum_{j=1}^qb_j\sigma_j^2)=\partial S_2-\Gamma_2;\,\,\partial(\sum_{k=1}^rc_k\sigma_k^3)=\Gamma_1-\Gamma_2.$$

Noticing that $\Phi$ is  well-defined on  $\sigma_i^1,\sigma_j^2,\sigma_k^3$, we can apply the Stokes theorem as follows:
\begin{equation*}
\begin{split}
0=\int_{\partial(\sum_{i=1}^pa_i\sigma_i^1)}d\Phi=&\int_{\partial S_1}\Phi-\int_{\Gamma_1}\Phi;\,\,\,\,\,
0=\int_{\partial(\sum_{j=1}^qb_j\sigma_j^2)}d\Phi=\int_{\partial S_2}\Phi-\int_{\Gamma_2}\Phi;\\
&0=\int_{\partial(\sum_{k=1}^rc_k\sigma_k^3)}d\Phi=\int_{\Gamma_1}\Phi-\int_{\Gamma_2}\Phi.
\end{split}
\end{equation*}
We conclude Lemma \ref{connect}. \,\,\,\,$\endpf$

\begin{definition}\label{tai}  Suppose $A=\big(X,W,(W_1,\cdots,W_l),(\gamma_1,\cdots,\gamma_m)\big)$ is a garden.
Let  $\Phi$ be a closed meromorphic $1$-form on $X$ with singularities on $W$.    We call $(d_1,\cdots,d_l)$ the short period vector of $\Phi$ with respect to $A$ if
	$$d_j=\int_{\partial S_j}\Phi\,\,\,{\text{for}}\,\,i=1,\cdots,l.$$ 
Here $S_j$, as in Lemma $\ref{connect}$, is a small analytic disc intersecting $W_j$ transversally at a smooth point $a\in W_j$; $\partial S_j$ is the counterclockwise oriented circle boundary of $S_j.$
\end{definition}
\begin{remark}
By Lemma $\ref{connect}$, the short period vector of a closed meromorphic $1$-form is well-defined.
\end{remark}

\begin{remark}\label{ret} The residue divisor {\rm Res}$(\Phi)$ of the closed meromorphic $1$-form $\Phi$ satisfies the following relation:
	$${\rm Res}(\Phi)=\sum_{j=1}^ld_jW_j,$$
	where $d_j$ is the $j$-component of the short period vector of $\Phi$.
\end{remark}
\begin{remark}\label{ess} The long period vectors (resp. the short period vectors) have different matrix representations with respect to different gardens;   however, they are related by basis changes of the vector space $\mathbb C^m$ (resp. the vector space $\mathbb C^l$). 
\end{remark}
\begin{definition}\label{lake pair} Let $A=\big(X,W,(W_1,\cdots,$ $W_l),(\gamma_1,\cdots,\gamma_m)\big)$ be a garden of $X$.  We say that  a closed meromorphic $1$-form or a closed anti-meromorphic $1$-form belongs to  $A$ if it is defined on $X$ with singularities on $W$. 
	
	Let $\Phi$ (resp.  $\widehat\Phi$) be a meromorphic $1$-form (resp. a anti-meromorphic $1$-form) belonging to $G$;  denote  by $(b_1,\cdots,b_m)$ (resp. $(\hat b_1,\cdots,\hat b_m)$) the long period vector of  $\Phi$ (resp.  $\widehat\Phi$) and $(c_1,\cdots,c_l)$ (resp. $(\hat c_1,\cdots,\hat c_l)$ ) the short period vector of $\Phi$ (resp.  $\widehat\Phi$) with respect to $A$.  We call $(\Phi,\widehat \Phi)$ a pair belonging to $A$ if and only if $b_i=-\hat b_i$ for $i=1,\cdots,m$ and $c_j=-\hat c_j$ for $j=1,\cdots,l.$ In this case, we call $\widehat \Phi$ an anti-meromorphic conjugate of $\Phi$ and $\Phi$ an meromprhic conjugate of  $\widehat\Phi$. For convenience, sometimes, we call $\Phi$ a conjugate of $\widehat \Phi$, and  vice versa.
\end{definition}

\begin{proposition}\label{pair} Let $A=\big(X,W,(W_1,\cdots,$ $W_l),(\gamma_1,\cdots,\gamma_m)\big)$ be  a garden.  Let $\widehat \Phi$ be a closed meromorphic $1$-form belonging to $A$. Then the conjugate of  $\widehat\Phi$ is unique modulo the differential of a meromorphic function. To be more precise, if $(\Phi_1,\widehat{\Phi})$ and  $(\Phi_2,\widehat{\Phi})$ are two pairs belonging to $A$, then $\Phi_1-\Phi_2=df$ where $f$ is a meromorphic function on $X$ with singularities on $W.$ 

\end{proposition}
\noindent{\bf {Proof:}} By assumption, $\Phi_1-\Phi_2$ is a closed meromorphic whose long period vector and short period vector both vanish. Then by Lemma $\ref{uniqueness}$ in the following, we draw the conclusion. $\endpf$

\begin{lemma}[Uniqueness Lemma]\label{uniqueness} Suppose that $\Phi$ is a closed meromorphic $1$-form with singularities on $W$. Assume the long period vector of $\Phi$ and the short period vector of $\Phi$ are both zero. Then $\Phi$ is the differential of a meromorphic function on $X.$
\end{lemma}
\noindent{\bf {Proof:}} Fix a point $p\in X\backslash W$. For each point $x\in X\backslash W$, we can find a smooth curve $\gamma_x$ connecting $p$ and $x.$ Define function $f(x)$ as follows:
$$f(x):=\int_{\gamma_x}\Phi.$$ We shall prove that $f(x)$ is well-defined, or equivalently, the integral is independent of the curve $\gamma_x$ connecting $p$ and $x$. It suffices to prove that for any smooth Jordan $\Gamma\subset X\backslash W$ based at $p$,
$$\int_{\Gamma}\Phi=0.$$

Let smooth 1-cycles $\tau_1,\cdots,\tau_m$ be a basis of the singular homology $H_1(X,\mathbb C)$ such that $\tau_i$ is a cycle contained in $X\backslash W$ with base point $p$ for $i=1,\cdots,m$. Then $\Gamma$ is homologous to $\sum_{i=1}^ma_i\tau_i.$ By Lemma \ref{smoothing} in Appendix, we can find finitely many smooth $2$-simplexes $\{\sigma_j\}_{j=1}^J$ so that the following properties hold:
\begin{enumerate}
	\item $\sigma_j:\Delta_2\rightarrow X$ is smooth   for $j=1,\cdots,J.$ Here $\Delta_2:=\{(x,y)\subset\mathbb R^2|0\leq x\leq 1,0\leq y\leq 1,x+y\leq 1\}$ is the standard 2-simplex.
	\item $\partial\big(\sum_{j=1}^J b_j\sigma_j\big)=\Gamma-\sum_{k=1}^{m}a_k\tau_k$, where $0\neq b_j\in\mathbb C$ for $j=1,\cdots,J.$ (By abuse of notation, we identify the image of $\Gamma$ with its corresponding $1$-cycle; the same for $\tau_k,k=1,\cdots,m.$)
	\item There are finitely many 1-simplexes  $\{\tilde\tau_l\}_{l=1}^L$ and finitely many 0-simplexes $\{A_n\}_{n=1}^N$  such that 
	\begin{equation*}
	\begin{split}
	\partial\sigma_{j}=\sum_{l=1}^Lc_{jl}\tilde\tau_{l}\,\,{\text{with}}\,\,c_{jl}\in\{0,1,-1\},\,\,{\text{for}}\,\,j=1,\cdots J,l=1,\cdots L;\\
	\partial\tilde\tau_{l}=\sum_{n=1}^Kd_{ln}A_{n}\,\,{\text{with}}\,\,d_{ln}\in\{0,1,-1\},\,\,{\text{for}}\,\,l=1,\cdots L,n=1,\cdots N. 
	\end{split}
	\end{equation*}
	Notice that $\{\Gamma,\tau_1,\cdots,\tau_m\}\subset\{\tilde\tau_l\}_{l=1}^L$ and $p\in\{A_n\}_{n=1}^N.$
	\item The above $\{\sigma_j\},\{\tilde\tau_l\}$ and $\{A_n\}$ are transversal to $W$ in the sense that: $A_n\notin W$ for $n=1,\cdots,N$; $\tilde\tau_l\cap W=\varnothing$ for $l=1,\cdots,L$; $\sigma_j\cap {\text{Sing}}(W)=\varnothing$ and $\sigma_j$ intersects $W$ transversally for $j=1,\cdots,J$. 
\end{enumerate}

For  $\sigma_j,j=1,\cdots,J,$ denote by $a^j_1,\cdots,a^j_{i_j}\subset\Delta_2$ the intersection points of $\sigma_j$ with $W$. Take small circles $S^j_1,\cdots,S^j_{i_j}\subset\Delta_2$ around  $a^j_1,\cdots,a^j_{i_j}$, respectively.
Notice that $\sigma_j^{*}(\Phi)$ is a well-defined closed $1$-form   on $\widehat\Delta_{2}^j:=\Delta_2\backslash\big(\cup _{q=1}^{j_i}S^j_q\big),j=1,\cdots,J.$ Applying  the Stokes theorem, we have 
\begin{equation}
\begin{split}
0=\int_{\widehat\Delta_{2}^j}\sigma_j^{*}(d\Phi)=\int_{\partial\widehat\Delta_{2}^j}\sigma_j^{*}(\Phi)=\sum_{l=1}^L\pm c_{jl}\int_{[0,1]}\tilde\tau_l^*(\Phi)+\sum_{q=1}^{i_j}\pm\int_{S^j_q}\sigma_j^*(\Phi),\,\,j=1,\cdots,J.
\end{split}
\end{equation}
When $S^j_{i}$ is small enough, by Lemma $\ref{connect}$ and the assumption that the short period vector of $\Phi$ is zero, we get
\begin{equation}
\sum_{q=1}^{i_j}\pm\int_{S^j_q}\sigma_j^*(\Phi)=0,\,\,j=1,\cdots,J.
\end{equation}
Since $\partial\big(\sum_{j=1}^J b_j\sigma_j\big)=\Gamma-\sum_{k=1}^{m}a_k\tau_k$, we obtain
$$\int_{\Gamma}\Phi=\sum_{k=1}^ma_k\int_{\tau_k}\Phi=0.$$
The last identity is because the long period vector of $\Phi$ is zero.

Therefore, we proved that $f(x)$ is well-defined.  It is easy to see that $f$ is a meromorphic function on $X$ with singularities on $W.$   Hence, we complete the proof of Lemma $\ref{uniqueness}$. $\endpf$

\begin{remark} Lemma $\ref{uniqueness}$ holds for closed anti-meromorphic $1$-forms in a similar way, and hence Proposition $\ref{pair}$  also holds for conjugates of a closed meromorphic $1$-form. That is, if there are two pairs $(\Phi,\widehat{\Phi}_1)$ and $(\Phi,\widehat{\Phi}_2)$, then  $\widehat \Phi_1-\widehat \Phi_2=df$ where $f$ is an anti-meromorphic function on $X$ with singularities  on $W.$ 
\end{remark}

\subsection{The existence of anti-meromorphic conjugates}
In this subsection we will prove a existence theorem for anti-meromorphic  conjugates of a closed meromorphic $1$-form belonging to a garden. We start with the following lemma.
\begin{lemma}\label{ge} Let $G=(X,W)$ be a garden. Let $\widehat W$ be an effective, ample divisor with the support contained in $W$.  Suppose $\overrightarrow{b}\in\mathbb C^m$ and $D\in\bigoplus_{i=1}^l\mathbb C  {W_i}$ with $c_1(D)=0$. Then there is a closed meromorphic $1$-form in $\Phi^1(*W)$ with long period vector $\overrightarrow{b}$ and residue divisor $D$. 
\end{lemma}
\noindent{\bf {Proof:}} By Theorem $\ref{WK}$, there is a closed meromorphic $1$-form $\Phi\in\Phi^1(W)$ with residue divisor $D.$ In the same way as Proposition $\ref{refined}$, we can find $\Psi$, a closed meromorphic $1$-form of the second kind, such that $\Psi\in\Phi^1(*W)$ and $\Psi+\Phi$ has long period vector $\overrightarrow{b}$.  Therefore, we proved Lemma $\ref{ge}$. $\endpf$


\medskip

Next we will prove the following existence theorem.
\begin{theorem}\label{con} Let $A=\big(X,W,(W_1,\cdots,W_l),(\gamma_1,\cdots,\gamma_m)\big)$ be a garden.  Let  $\Phi$  be a closed meromorphic $1$-form belonging to  $A$. There exists a conjugate of $\Phi$ belonging to  $A$.
\end{theorem}
\noindent{\bf {Proof :}} Denote by $(b_!,\cdots,b_m)$ the long period vector of $\Phi$. Let  Res$(\Phi)$ be the residue divisor of $\Phi$ with the form  Res$(\Phi)=\sum_{i=1}^la_iW_i$. Since $\Phi$ is a closed meromorphic $1$-form belonging to $A$, the first Chern class of Res$(\Phi)$ is zero, that is, $\sum_{i=1}^la_i\cdot c_1(W_i)=0\in H^2(X,\mathbb C)$.   Noticing that $\sum_{i=1}^l(-\overline{a_i})\cdot c_1(W_i)=0\in H^2(X,\mathbb C)$,  by Lemma $\ref{ge}$, we can derive a closed meromorphic $1$-form $\Psi$ belonging to $G$, with residue divisor $\sum_{i=1}^l(-\overline{a_i})\cdot W_i$ and long period vector $(-\overline{b_1},\cdots,-\overline {b_m})$. 

Take the complex conjugate $\overline{\Psi}$ of $\Psi$; that is, if we write $\Psi=\psi_1+\sqrt{-1}\psi_2$ where $\psi_1$ and $\psi_2$ are closed real $1$-forms with singularities on $W$, then $\overline{\Psi}=\psi_1-\sqrt{-1}\psi_2$. It is easy to verify that $\overline{\Psi}$ is an anti-meromorphic conjugate of $\Phi$.  We complete the proof of Theorem \ref{con}. \,\,\,\,  $\endpf$
\medskip

\subsection{Constructing pluriharmonic functions with log poles}
In this subsection we will construct pluriharmonic functions with log poles by integrating the sum of pairs.
\begin{theorem}\label{epluri} Let $A=\big(X,W,(W_1,\cdots,W_l),(\gamma_1,\cdots,\gamma_m)\big)$ be a garden. Assume $(\Phi,\widehat \Phi)$ is a pair belonging to $G$. Then the following integral is well-defined on $X\backslash W$ and gives a pluriharmonic function with singularities on $W:$
\begin{equation}
h_{\Phi\widehat{\Phi}}(z):=\int_{\gamma^p_{z}}(\Phi+\widehat{\Phi}).
\end{equation}
Here $p\in X\backslash W$ is a fixed point and $\gamma_{z}^p$ is a smooth curve connecting $p$ and $z\in X\backslash W.$ 
\end{theorem}
\noindent{\bf {Proof:}} It suffices to prove that along each Lipschitz loop $\gamma\subset X\backslash W$ the following holds:
\begin{equation}
\int_{\gamma}(\Phi+\widehat{\Phi})=0\,\,.
\end{equation}
In the same way as the proof of Lemma $\ref{uniqueness}$, we can find complex numbers $\{a_i\}_{i=1}^m$ and $\{b_j\}_{j=1}^l$ so that 
\begin{equation}\label{80}
\int_{\gamma}(\Phi+\widehat{\Phi})=\sum_{i=1}^ma_i\int_{\gamma_i}(\Phi+\widehat{\Phi})+\sum_{j=1}^lb_j\int_{\partial S_j}(\Phi+\widehat{\Phi}).
\end{equation}
Here $\partial S_j$ is the circle boundary of a small analytic disc intersecting $W_j$ at its smooth point.
Notice that the sum of the long period vector (resp. the short period vector)  of $\Phi$ and the long period vector (resp. the short period vector)  of $\widehat \Phi$ is zero.
Then each term on the right hand side of formula $(\ref{80})$ vanishes. 
We draw the conclusion.  $\endpf$

\begin{definition}
	We call $h_{\Phi\widehat{\Phi}}$ a pluriharmonic functions coming from the pair $(\Phi,\widehat \Phi)$.
\end{definition}
\noindent{\bf Proof of Theorem \ref{iepluri} :} Take a garden $A=\big(X,W,(W_1,\cdots,W_l),(\gamma_1,\cdots,\gamma_m)\big)$ of $X$. By Theorem \ref{con}, for each closed meromorphic $1$-form $\Phi$ with poles on $W$, there exists a conjugate $\widehat{\Phi}$ of $\Phi$. By Theorem \ref{epluri}, we conclude Theorem \ref{iepluri}.\,\,\,$\endpf$
\medskip

Let $G=(X,W)$ be a garden.
Denote by $M_G$ the vector space of meromorphic functions on $X$ with singularities on $W$; denote by $\widehat M_G$ the vector space of anti-meromorphic functions with singularities on $W$. Define $T_G:=M_G+\widehat M_G.$
Next we will establish the uniqueness of pluriharmonic functions and estimate the dimension of the vector spaces of pluriharmonic functions coming from pair (modulo $T_G$).
\begin{theorem}\label{upluri} Let $G=(X,W)$ be a garden. Assume $(\Phi_1,\widehat \Phi_1)$ and $(\Phi_2,\widehat \Phi_2)$ are two pairs belonging to $G$ . Then $h_{\Phi_1\widehat{\Phi}_1}=h_{\Phi_2\widehat{\Phi}_2}$ (modulo $T_G$) if and only if  $\Phi_1$ and $\Phi_2$ have the same long period vector and the same short period vector.
\end{theorem}
\noindent{\bf {Proof:}} Assume that $\Phi_1$ and $\Phi_2$ have the same long period vector and the same short period vector.  Then so do $\widehat\Phi_1$ and $\widehat\Phi_2$. By Lemma $\ref{uniqueness}$, the following functions are well-defined on $X:$
$$f:=\int_{\gamma_z^p}(\Phi_1-\Phi_2);\,\,\,\,\,\,\,\bar h:=\int_{\gamma_z^p}(\widehat\Phi_1-\widehat\Phi_2).$$
Noticing that $f$ is meromorphic and $\bar h$ is anti-meromorphic, we conclude that $h_{\Phi_1\widehat{\Phi}_1}=h_{\Phi_2\widehat{\Phi}_2}$ (modulo $T_G$).

On the other hand if $h_{\Phi_1\widehat{\Phi}_1}=h_{\Phi_2\widehat{\Phi}_2}$ (modulo $T_G$), then there are meromorphic function $f$ and anti-meromorphic function $\bar h$ on $X$ such that 
\begin{equation}\label{81}
\begin{split}
\int_{\gamma^p_{z}}(\Phi_1+&\widehat{\Phi}_1)-\int_{\gamma^p_{z}}(\Phi_2+\widehat{\Phi}_2)=\overline h-f\,\,{\text{or equivalently}}\\
&f+\int_{\gamma^p_{z}}(\Phi_1-\Phi_2)=\overline h+\int_{\gamma^p_{z}}(\widehat{\Phi}_2-\widehat{\Phi}_1).
\end{split}
\end{equation}
Notice that the left hand side of formula $(\ref{81})$ is holomorphic and the right hand side is anti holomorphic,  away from the singularities of $\Phi_1,\widehat{\Phi}_1,\Phi_2,\widehat{\Phi}_2,f$ and $\overline h.$ Then both sides are locally constants. Taking differential of the left hand side, we have the following equality on $X$:
\begin{equation}
df+\Phi_1-\Phi_2=0.
\end{equation}
Integrating $df+\Phi_1-\Phi_2$ along $\gamma_i$ and $\partial S_j$, we conclude that the the long period vector (resp. the short period vector) of $\Phi_1$ and  that of $\Phi_2$ are the same. We complete the proof. $\endpf$

\begin{theorem}\label{dimension}Let $G=(X,W)$  be a garden. Denote by $k$ the dimension of the kernel of the double delta map $\delta^1\circ\Delta^0:H^0(X,R(W))\rightarrow H^2(X,\mathbb C).$  Then the dimension of the vector space of the pluriharmonic functions coming from pairs with singularities on $W$ (modulo $T_G$) is $k+\dim H^1(X,\mathbb C)$.
\end{theorem}
\noindent{\bf {Proof :}} By Theorem $\ref{cloister}$, the dimension of the closed meromorphic $1$-forms of second kind modulo exact form is $\dim H^1(X,\mathbb C)$. By Theorem  $\ref{HA}$, the dimension of the closed meromorphic $1$-forms modulo the closed meromorphic $1$-forms of the second kind is $k.$  By Theorem  $\ref{epluri}$, we can construct pluriharmonic functions for every closed meromorphic $1$-form; by Theorem $\ref{upluri}$, these functions are linear independent modulo $T_G.$  Therefore we conclude Theorem \ref{dimension}.
$\endpf$

\begin{remark} Denote by $V$ be the vector space of pluriharmonic functions coming from a closed meromorphic $1$-form. Then $\dim V/T_G=\infty$. On the other hand, denote by $V_2$ the vector space of pluriharmonic functions coming from  a closed meromorphic $1$-form of the second kind.  Then $\dim V_2/T_G=\dim H^1(X,\mathbb C)<\infty$.
\end{remark}

\noindent{\bf {Proof of Theorem $\ref{121}$ :}}
Let $h$ be a pluriharmonic function on $X$ of local form $(\ref{log})$.  Then  $h$ locally takes the form of
	\begin{equation}\label{q2}
	h(z)={g_1(z)}+{g_2(\bar z)}+\sum_{i=1}^la_i\log|f_i|^2\,\,\,,
	\end{equation}
where $a_1,\cdots,a_l$ are constants and $g_1,g_2,f_1,\cdots,f_l$ are meromorphic functions.  Taking differential of $h$, we have that locally
\begin{equation}
\begin{split}
&dh=\partial{g_1(z)}+\sum_{i=1}^la_i\frac{\partial f_i}{f_i}+\bar{\partial}{g_2(\bar z)}+\sum_{i=1}^la_i\frac{\bar\partial \bar{f_i}}{\bar f_i}\,\,\,;\\
&\Phi:=\partial{g_1(z)}+\sum_{i=1}^la_i\frac{\partial f_i}{f_i};\,\,\,\widehat{\Phi}:=\bar{\partial}{g_2(\bar z)}+\sum_{i=1}^la_i\frac{\bar\partial \bar{f_i}}{\bar f_i}.
\end{split}
\end{equation}
Notice that $\Phi$ is a closed meromorphic $1$-form and  $\widehat\Phi$ is a closed anti-meromorphic $1$-form. Moreover, it is clear that $dh=\Phi+\widehat{\Phi}$ holds globally on $X$; that is, the definition of $\Phi$ and $\widehat{\Phi}$ does not depend on the choice of the local charts.
If there is no log term in formula $(\ref{q2})$, then  $\Phi$ is a closed meromorphic $1$-form of the second kind.  One can show that $h$ is the integral of the sume of $\Phi$ and $\widehat\Phi$ in the same way as Theorem \ref{epluri}.

Define $\kappa$ by $\kappa(h)=\Phi$ for every $h\in Ph(X)$. It is clear that $\kappa$ is well-defined. Next we will show that $\kappa$ is injective; that is, for $h_1,h_2\in Ph(X)$ if $\kappa(h_1)=\kappa(h_2)$, then $h_1-h_2\in (K(X)+\overline K(X))$.  Notice that the singularities of $h_i$ is contained in a divisor of $X$ for $i=1,2$. Take an ample, reduced divisor $W$ of $X$ such that the singularities of $h_i$ is contained in $W$ for $i=1,2$; hence $(X,W)$ is a garden. Then by Theorem \ref{upluri}
, we conclude the injectivity of $\kappa$.

Finally, we will show that $\kappa$ is surjective. Let $\Phi\in  H^0(X,\Phi^1(*)).$  Take an ample, effective, reduced divisor $W$ of $X$ such that the singularities of $h_i$ is contained in $W$ for $i=1,2$; hence $(X,W)$ is a garden. Then by Theorem \ref{con}, we can find a conjugate of $\Phi$. By Theorem \ref{epluri}, we conclude the surjectivity.

Similarly, we can show that $\kappa_0$ is an isomorphism.  Hence we complete the proof of Theorem \ref{121}.\,\,$\endpf$

\begin{remark}[\cite{Y}] The log term in the above formula has an interesting explanation related to superfluid vortices in physics. A single vortex has energy growing asymptotically like $\log z_i$ and hence unstable. In low temperature, there will be no free vortices, only clusters of zero total vorticity.
\end{remark}

\appendix
\section{Appendix I: Good Cover Lemma}

In this appendix we will present a detailed proof of the following well known lemma for reader's convenience. 

\begin{lemma}[Good cover lemma]\label{GC} Let $X$ be a compact complex manifold $X$. Let $\mathcal V$ be an open cover of $X$. Then, there is a finite open cover $\mathcal U:=\{U_i\}_{i=1}^M$ of $X$ such that the following properties hold:
	\begin{enumerate}
		\item $U_i$ is a Stein space for $i=1,\cdots,M$;
		\item Each intersection of $U_i$ is contractible  for $i=1,\cdots,M$, if it is nonempty;
		\item $\mathcal U$ is a refinement of $\mathcal V$.
	\end{enumerate}
\end{lemma}

In order to prove Lemma \ref{GC}, we will first prove the following lemma.
\begin{lemma}\label{convex}
Let $X$ be a complex manifold of complex dimension $m$. Suppose there are two holomorphic local parametrizations of $X$ as 
\begin{equation}
\Phi_i:V_i\rightarrow U_i\subset X,\,\,i=1,2,
\end{equation}
where $V_i$ is an open set in $\mathbb C^m$ and $U_i$ is biholomorphic to $V_i$. Moreover, for each point $x=(x_1,\cdots,x_m)\in V_i$ and  real number $r>0$, denote by $B_{i,r}(x)$ the following complex ball in $V_i$: 
\begin{equation}\label{ball}
B_{i,r}(x):=\{(y_1,\cdots,y_m)\in V_i\big||y_1-x_1|^2+|y_2-x_2|^2+\cdots |y_m-x_m|^2<r^2\}.
\end{equation}
Then for each point $x\in \Phi_1^{-1}(U_1\cap U_2)$, there is a positive number $R$ such that  for any positive  number $r<R$ the complex ball $B_{1,r}(x)\subset\Phi_1^{-1}(U_1\cap U_2)$ and, moreover,  its biholomorhpic image  $(\Phi_2^{-1}\circ\Phi_1)(B_{1,r}(x))$ is a real convex set in $V_2$.
\end{lemma}
\noindent{\bf {Proof :}} Fix a point $x\in \Phi_1^{-1}(U_1\cap U_2)$. Denote by $(z_1,\cdots,z_m)$ the complex coordinates for $V_1$; denote by $(w_1,\cdots,w_m)$ the complex coordinates for $V_2$.  Without loss of genearlity, we can assume that the point $x\in \Phi_1^{-1}(U_1\cap U_2)$ has complex coordinates $(0,\cdots,0)$ in $V_1$ and the point $y:=(\Phi_2^{-1}\circ\Phi_1)(x)$ has  complex coordinates $(0,\cdots,0)$ in $V_2$. For convenience,  denote by $\phi$ the restriction of the holomorphic map $\Phi_1^{-1}\circ\Phi_2$ to $\Phi_2^{-1}(U_1\cap U_2)$, that is, 
\begin{equation}
\begin{split}
\phi:\Phi_2^{-1}(U_1\cap U_2)&\rightarrow \Phi_1^{-1}(U_1\cap U_2),\\
(w_1,\cdots,w_m)&\mapsto(z_1,\cdots,z_m)=(\Phi_1^{-1}\circ\Phi_2)(w_1,\cdots,w_m).
\end{split}
\end{equation}
Note that $\phi(y)=x.$
Restricted to a small ball $B_{2,s}(y)\subset\Phi_2^{-1}(U_1\cap U_2)$, $s>0$, we have the following Taylor expansion:
\begin{equation}
z_1=\sum_{i=1}^na_{1i}w_i+O(|w|^2);
z_2=\sum_{i=1}^na_{2i}w_i+O(|w|^2);
\cdots;
z_m=\sum_{i=1}^na_{mi}w_i+O(|w|^2).
\end{equation}
Since $\phi$ is a biholomorphic map,  the matrix $A:=(a_{ij})_{i,j=1}^{m}$ is non-degenerate. Define a real analytic function $F(w_1,\cdots,w_m)$ in $B_{2,s}(y)$ as follows: 
\begin{equation}
\begin{split}
F(w_1,&\cdots,w_m)=\sum_{i=1}^m|z_i(w_1,\cdots,w_m)|^2=\sum_{i=1}^m\big|\sum_{j=1}^na_{ij}w_j+O(|w|^2)\big|^2\\
=&\sum_{i=1}^m\big|\sum_{j=1}^na_{ij}w_j\big|^2+O(|w|^3)=(\overline w_1,\cdots,\overline w_m)\cdot\overline {A}^T\cdot A\cdot(w_1,\cdots,w_m)^T+O(|w|^3).
\end{split}
\end{equation}
Take a complex ball $B_{1,\widetilde R}(x)$ centered at $x$ such that 
$B_{1,\widetilde R}(x)\subset \phi(B_{2,s}(y));$  then, the set $\phi^{-1}(B_{1,r}(x))$ in $V_2$  is the same as the set $\{F<r^2\}$ for $0<r<\widetilde R$.

Since ${A}^T\cdot A$ is a strictly positive definite Hermitian matrix, we can find an $m\times m$ unitary matrix $U$ such that 
\begin{equation}
{U}^T\cdot{A}^T\cdot A\cdot \overline U=\left(\begin{array}{cccc} a_1 & 0 &\cdots & 0\\0 & a_2 &\cdots & 0\\\cdots&\cdots&\cdots&\cdots\\ 0 &\cdots& 0& a_m \end{array}\right)=:D, 
\end{equation}
where $a_i$ is a positive real number for $i=1,\cdots,m.$  Take an unitary change of coordinates for $B_{1,\widetilde R}(x)$ as
\begin{equation}
(\widetilde w_1,\cdots,\widetilde w_m)=(w_1,\cdots,w_m)\cdot U.
\end{equation}
Then $F$ in the new coordinates takes the form of 
\begin{equation}
\sum_{i=1}^ma_i|\widetilde w_i|^2+O(|\widetilde w|^3);
\end{equation}
we denote it by $\widetilde F$.

Notice that a unitrary transformation of $\mathbb C^m$ does not change the shape of  domains. Therefore, in order to prove Lemma \ref{convex}, it suffices to prove that there exists a positive number $R$ such that $\{\widetilde F<r^2\}$ is a real convex domain in $V_2$ for each $0<r<R$.  Denote the real coordinates of $(\widetilde w_1,\cdots,\widetilde w_m)$ by $(x_1,\cdots,x_m,y_1,\cdots,y_m)$ where $\widetilde w_i=x_i+\sqrt{-1}y_i$ for $i=1,\cdots,m.$ Computation of the real Hessian of $\widetilde F$ yields that
\begin{equation}
{\rm Hess}(\widetilde F)=\left(\begin{array}{cc} (\frac{\partial^2\widetilde F}{\partial x_i\partial x_j}) & (\frac{\partial^2\widetilde F}{\partial x_i\partial y_j})\\ (\frac{\partial^2\widetilde F}{\partial y_i\partial x_j}) &  (\frac{\partial^2\widetilde F}{\partial y_i\partial y_j})\\ \end{array}\right)=
\left(\begin{array}{cc} D & 0 \\0&D \\ \end{array}\right)+O(|\widetilde w|);
\end{equation}
hence the real Hessian of $\widetilde F$ is strictly positive in a small neighborhood of $y$. Therefore, we can choose a positive number $R$ such that $\{\widetilde F<r^2\}$ is a real convex domain in $V_2$ for each $0<r<R$.\,\,\,\,\,\,\,$\endpf$

Now we turn to the proof of Lemma \ref{GC}.
\medskip

\noindent{\bf {Proof of Lemma \ref{GC}:}} Fix a Riemannian metric $g$ on $X$. For any points $x,y\in X$, denote by $d_g(x,y)$ the distance between them with respect to $g$. Since $X$ is compact, we can have finitely many local parametrizations $\{\phi_i:V_i\rightarrow U_i\}_{i=1}^M$ such that the following properties hold.
\begin{enumerate}
	\item $V_i=\{z=(z_1,\cdots,z_m)\in\mathbb C^m\big|d_g(\phi_i(0),\phi_i(z))<4\}$ for $i=1,\cdots,M$.
	\item $U_i\subset X$ and $\phi_i$ is the biholomorphic map between $V_i$ and $U_i$ for $i=1,\cdots,M$.
	\item For $i=1,\cdots,M$, denote by $V_{i,1}$, $V_{i,2}$, $U_{i,1}$ and $U_{i,2}$ the set $\{z\in\mathbb C^m\big|d_g(\phi_i(0),\phi_i(z))<1\}$, $\{z\in\mathbb C^m\big|d_g(\phi_i(0),\phi_i(z))<2\}$, $\phi_i(V_{i,1})$ and $\phi_i(V_{i,2})$, respectively. Then, $\{U_{i,1}\}_{i=1}^M$ is an open cover of $X$.
\end{enumerate}
For convenience, denote by $I$ the index set $\{1,2,\cdots,M\}$. Morevover, similar to formula $(\ref{ball})$, for $i=1,\cdots, M$, $0<r<1$ and each point $x\in V_{i,1}$ with complex coordinates $(x_1,\cdots,x_m)$, we denote by $B_{i,r}(x)$ the set $\{(y_1,\cdots,y_m)\in V_i\big||y_1-x_1|^2+|y_2-x_2|^2+\cdots |y_m-x_m|^2<r^2\}$.

Next, we will construct a special open cover of $V_{i,1}$ for $i=1,\cdots,M$. Fix $i$ and, for each point $x\in V_{i,1}$, let $I^i_x$ be the index set define by
\begin{equation}
I_x:=\{j\in I\big|\phi_i(x)\in U_{j,2}\}.
\end{equation}
Then, for each $j\in I^i_x$, we can apply Lemma \ref{convex} and derive a positive number $R_j$ such that for any $0<r<R_j$ the following properties hold:
\begin{enumerate}[\label={}]
	\item $\bullet$ $\phi_i(B_{i,r}(x))\subset U_{j,2}$;
	\item $\bullet$ $\phi_j^{-1}(\phi_i(B_{i,r}(x)))$ is a real convex set in $V_{j,2}$.
\end{enumerate}
Choose a positiv number $r_{x}$ such that $r_{x}<\min_{j\in I^i_x}\{R_j\}$, the diameter of $B_{i,r_x}(x)$ with respect to $d_g$ is less than $\frac{1}{3}$, and $\phi_i(B_{i,r}(x)$ is contained in a certain open set in open cover $\mathcal V$. Notice that $B_{i,r_{x}}(x)$ is an open neighborhood of $x$; moreover,  the following properties hold,
\begin{enumerate}[\label={}]
	\item $	(\ast) $ $\phi_i(B_{i,r_{x}}(x))\subset U_{j,2}$ for each $j\in I^i_x$,
	\item $(\star) $ $\phi_j^{-1}(\phi_i(B_{i,r_{x}}(x)))$ is a real convex set in $V_{j,2}$ for each $j\in I^i_x$.
\end{enumerate}
It is clear that $\big\{B_{i,r_{x}}(x),x\in V_{i,1}\big\}$ is an open cover of $V_{i,1}$.

Since $\bigcup_{i=1}^M\big\{\phi_i(B_{i,r_{x}}(x)),x\in V_{i,1}\big\}$ is an open cover of $X$,  we can find a finite subcover $\Omega$ of $X$ with the form of
\begin{equation}
\Omega=\bigcup_{i=1}^M\big\{\phi_i(B_{i,r_{a_{ij}}}(a_{ij})),a_{ij}\in V_{i,1},j=1,\cdots,N_i\big\}.
\end{equation}

We claim that the cover $\Omega$ satiesfies the properties required in Lemma \ref{GC}. The first and the third properties are clear, for each open set in $\Omega$ is biholomoprhic to a complex ball in $\mathbb C^m$ and is contained in a certain open set in $\mathcal V$. Since convext sets are contractible, in order to establish the second property, it suffices to show that each nonempty interesections of open sets in $\Omega$ is boholomorphic to a real convex set in $\mathbb C^m.$  

Let $B$ be an nonempty intersection of open sets in $\Omega$. Without loss of generality, by rearranging the indices, we can assume that $B$ takes the form of 
\begin{equation}
B=\bigcap_{i=1}^M\bigcap_{j=1}^{l_i}\phi_i(B_{i,r_{a_{ij}}}(a_{ij})),
\end{equation}
where $l_i\leq N_i$. Without loss of generality, we can further assume $l_1\geq 1$. Since $B$ is nonempty, by the construction we have $a_{ij}\in V_{1,2}$. Therefore, $1\in I_{a_{ij}}$ for $i=1,\cdots,M$ and $j=1,\cdots,l_i.$ By property $(\ast)$ and property $	(\star)$,  we have that $\phi_i(B_{i,r_{x}}(x))\subset U_{1,2}$  and $\phi_1^{-1}(\phi_i(B_{i,r_{a_{ij}}}(a_{ij})))$ is a real convex set in $V_{1,2}$ for each  for $i=1,\cdots,M$ and $j=1,\cdots,l_i$. Since the intersection of convex sets are convex, $B$ is a convex set.  

Therefore, we complete the proof of Lemma \ref{GC}.\,\,\,\,\,\,\,$\endpf$
\medskip

We also include a proof of the following very good cover lemma for the completeness.
\begin{lemma}[Very good cover]\label{VGC} Let $X$ be a compact complex manifold and $W$ be a normal crossing divisor on $X$. Let $\mathcal V$ be an open cover of $X$. Then, there is a finite open cover $\mathcal U:=\{U_i\}_{i=1}^M$ of $X$ such that the following properties hold:
	\begin{enumerate}
		\item $U_i$ is a Stein space for $i=1,\cdots,M$;
		\item Each intersection of $U_i$, $i=1,\cdots,M,$  is contractible, if it is nonempty;
		\item There are only finitely many irreducible components of $U_i\bigcap W$ for $i=1,\cdots,M$; moreover, if $W_{ij}$ is an irreducible component of $U_i\cap W$,  then $W_{ij}$ is contractible;
		\item $\mathcal U$ is a refinement of $\mathcal V$.
	\end{enumerate}
\end{lemma}
\noindent{\bf Proof of Lemma \ref{VGC} :}
We first prove the following Claim.
\medskip

{\bf {Claim :}} Let $U$ be an open set in $\mathbb C^m$. Suppose that $W\subset U$ is a reduced analytic subvariety  of codimension one  with normal crossing singularities  only.  For each point $x\in U$, there is a Euclidean complex ball $B_{r}(x)\subset U$ such that there are only finitely many irreducible components in the irreducible decomposition of $W\bigcap B_{r}(x)$ and each irreducible component is contractible.
\medskip

{\bf {Proof of Claim :}} 
Since the case $x\in U\backslash W$ is trivial, let $x\in W\bigcap U$; without loss of generality, we assume $x=(0,\dots,0)\in U$. By taking the complex ball $B_R:=\{z\big||z|^2<R^2\}$ with radius $R>0$ small enough, we can assume that $B_R\subset U$ and $W$ is defined in $B_R$ by equation $\prod_{i=1}^lf_i=0$, where $f_1,\cdots,f_l$ are  holomorphic functions in $B_R$  such that they are irreducible, pairwise coprime and vanishing at $(0,\cdots,0)$. Without loss of generality, we assume that $f_1,\cdots,f_l$  take the  form of
\begin{equation}\label{Wi}
f_i(z)=\sum_{j=1}^ma_{ij}z_j+O(|z|^2)
\end{equation}
with $a_{ij}\in\mathbb C$ for $i=1\cdots,l,\,\,j=1,\cdots,m$ and $ z\in B_R$.
Since $W$ is a normal crossing divisor, $f_i$ has a nonzero linear part for $i=1,\cdots,l.$

Define $W_i:=\{f_i=0\}\bigcap B_R$. It is clear that $\bigcup_{i=1}^l W_i=W\bigcap B_R.$ Hence it suffices to prove that each $W_i$ is contractible when $R$ is small enough for $i=1,\cdots,l$. In the following, we will prove this for $W_1;$ the proof for others is the same which we will omit. 

Without loss of generality, we assume that $a_{11}=-1$ in formula (\ref{Wi}). By the implicit function theorem, we can solve $z_1$ in terms of $(z_2,\cdots,z_m)$ as 
\begin{equation}
z_1=h(z_2,\cdots,z_m)=a_{12}z_2+a_{13}z_3+\cdots+a_{1m}z_m+O(|z_2|^2+\cdots+|z_m|^2).
\end{equation}
Shrinking the radius $R$, we have a parametrization $(G,\widetilde W)$ of $W_1$ as 
\begin{equation}
\begin{split}
G:\,\,\,\widetilde W&\rightarrow W_1\subset\mathbb C^m,\\
(z_2,\cdots,z_m)&\mapsto(h(z_2,\cdots,z_m),z_2,z_3,\cdots,z_m),
\end{split}
\end{equation}
where $\widetilde W$ is a domain in $\mathbb C^{m-1}$ defined by 
\begin{equation}
{\rho}(z_2,\cdots,z_m):=\big|h(z_2,\cdots,z_m)\big|^2+\sum_{i=2}^m\big|z_i\big|^2<R^2.
\end{equation} Since $W_1$ and $\widetilde W$ are biholomorphic, it suffices to show that $\widetilde W$ is contractible. Computation yields that
\begin{equation}
\begin{split}
{\rho}&=\big|a_{12}z_2+a_{13}z_3+\cdots+a_{1m}z_m+O(|z|^2)\big|^2+|z_2|^2+\cdots+|z_m|^2\\
&=|a_{12}z_2+a_{13}z_3+\cdots+a_{1m}z_m|^2+|z_2|^2+\cdots+|z_m|^2+O(|z|^3).
\end{split}
\end{equation}
Taking derivatives, we get
\begin{equation}
\begin{array}{lll}
&\frac{\partial{\rho}}{\partial z_2}&=(a_{12}\bar a_{12}+1)\bar z_2+a_{12}\bar a_{13}\bar z_3+\cdots+a_{12}\bar a_{1m}\bar z_{m}+O(|z|^2),\\
&\frac{\partial{\rho}}{\partial z_3}&=a_{13}\bar a_{12}\bar z_2+(a_{13}\bar a_{13}+1)\bar z_3+\cdots+a_{13}\bar a_{1m}\bar z_{m}+O(|z|^2),\\
&&\cdots\\
&\frac{\partial{\rho}}{\partial z_m}&=a_{1m}\bar a_{12}\bar z_2+a_{1m}\bar a_{13}\bar z_3+\cdots+(a_{1m}\bar a_{1m}+1)\bar z_{m}+O(|z|^2).\\
\end{array}
\end{equation}
Multiplying $\frac{\partial{\rho}}{\partial z_i}$ by $z_i$, $i=1,\cdots,m$, and summing up the products, we obtain 
\begin{equation}
\begin{split}
\sum_{i=1}^m\frac{\partial{\rho}}{\partial z_i}z_i=|\sum_{i=1}^ma_{1i}z^i|^2+\sum_{i=1}^m|z_i|^2+O(|z|^3).
\end{split}
\end{equation}
Then for $R>0$ small enough
$$|\sum_{i=1}^m\frac{\partial{\rho}}{\partial z_i}z_i|>\frac{1}{2}\sum_{i=2}^m|z_i|^2\,\,\,{\text{for}}\,\,(z_2,\cdots,z_m)\in\widetilde W;$$
hence $\nabla{\rho}$ is nonzero except at the origin.  By Morse theory $\widetilde W$ is contractible.  We complete the proof of the claim. $\endpf$
\medskip

The remaining of the proof is exact the same as the proof of Lemma \ref{GC} except that when choosing $B_{i,r_x}(x)$, we require the following condition in addition
\begin{equation}
\phi_i(B_{i,r_x}(x))\bigcap W=\bigcup_{j=1}^{l_i}W_{j}
\end{equation}
where $W_j$ is irreducible and contractible for $j=1,\cdots,l_i<\infty$.\,\,\,\,$\endpf$

\section{Appendix II: Existence of a smooth, transversal two-chain}

In this appendix we give a detailed proof of the following lemma:
\begin{lemma}\label{smoothing}
	Let $X$ be a compact algebraic manifold, $p\in X$ a fixed base point and $W$ a  normal crossing divisor. Suppose $\tilde\tau_1,\cdots,\tilde\tau_m$ is a basis of the singular homology $H_1(X,\mathbb C)$, where
	\begin{displaymath}
	\tilde\tau_k:[0,1]\rightarrow X\,\, {\text{is smooth with}} \,\, \tilde\tau_k(0)=\tilde\tau_k(1)=p \,\,{\text{for}}\,\, k=1,\cdots,m.
	\end{displaymath}
	Assume $\tilde\tau$ is also a one-cycle, where $\tilde\tau:[0,1]\rightarrow X$ is smooth with $\tilde\tau(0)=\tilde\tau(1)=p$. Moreover we assume that  $\tilde\tau([0,1])\cap W=\varnothing$ and $\tilde\tau_k([0,1])\cap W=\varnothing$ for $k=1,\cdots,m.$
    Then if $\tilde\tau$ is homologous to $\sum_{k=1}^{m}a_k\tilde\tau_k$ with constants $\{a_k\}_{k=1}^m\subset\mathbb C$, we can find finitely many smooth  2-simplexes $\{\sigma_j\}_{j=1}^J$ such that the following properties are satisfied:
    \begin{enumerate}
    	\item $\sigma_j:\Delta_2\rightarrow X$ is smooth   for $j=1,\cdots,J.$ Here $\Delta_2:=\{(x,y)\subset\mathbb R^2|0\leq x\leq 1,0\leq y\leq 1,x+y\leq 1\}$ is the standard 2-simplex.
    	\item $\partial\big(\sum_{j=1}^J b_j\sigma_j\big)=\tilde\tau-\sum_{k=1}^{m}a_k\tilde\tau_k$, where $0\neq b_j\in\mathbb C$ for $j=1,\cdots,J.$ 
    	\item There are finitely many 1-simplexes $\{\tau_l\}_{l=1}^L$ and finitely many 0-simplexes $\{A_n\}_{n=1}^N$  such that 
    	\begin{equation*}
    	\begin{split}
    	\partial\sigma_{j}=\sum_{l=1}^Lc_{jl}\tau_{l}\,\,{\text{with}}\,\,c_{jl}\in\{0,1,-1\}\,\,{\text{for}}\,\,j=1,\cdots J,l=1,\cdots L;\\
    	\partial\tau_{l}=\sum_{n=1}^Kd_{ln}A_{n}\,\,{\text{with}}\,\,d_{ln}\in\{0,1,-1\}\,\,{\text{for}}\,\,l=1,\cdots L,n=1,\cdots N. 
    	\end{split}
    	\end{equation*} 
    	Notice that $\{\tilde{\tau},\tilde\tau_1,\cdots,\tilde{\tau}_m\}\subset\{\tau_j\}_{l=1}^L$ and $p\in\{A_n\}_{n=1}^N$.
    	\item The above $\{\sigma_j\},\{\tau_l\}$ and $\{A_n\}$ are transversal to $W$ in the sense that: $A_n\notin W$ for $n=1,\cdots,N$  and $\tau_l\cap W=\varnothing$ for $l=1,\cdots,L$; $\sigma_j\cap {\text{Sing}}(W)=\varnothing$ and $\sigma_j$ intersects $W$ transversally for $j=1,\cdots,J$. 
    \end{enumerate}
\end{lemma}
{\bf {Proof:}} Since $\tilde\tau$ is homologous to $\sum_{k=1}^{m}a_k\tilde\tau_k$, there are finitely many $2$-simplexes $\{\sigma_j\}_{j=1}^J$ such that $\partial\big(\sum_{j=1}^J b_j\sigma_j\big)=\tilde\tau-\sum_{k=1}^{m}a_k\tilde\tau_k$, where $0\neq b_j\in\mathbb C$ for $j=1,\cdots,J.$

The basic idea of the proof is to perturb the simplexes homotopically, thicken the lower dimensional simplexes and extend the perturbation from simplexes with lower dimension  to simplexes with higher dimension. Notice that throughout the proof the $0$-simplex $p$ and the $1$-simplexes $\{\tilde{\tau},\tilde\tau_1,\cdots,\tilde{\tau}_m\}$ are kept the same without any change. We divide the proof into seven steps as follows.

Step 1: Perturb $\{A_n\}$ so that $A_n\notin W$ for $n=1,\cdots,N$. To be more precise, we shall construct homotopies $I_n:[0,1]\rightarrow X,n=1,\cdots,N$ such that:
\begin{enumerate}
	\item  $I_n$ is smooth;
	\item  $I_n(0)=A_n,I_n(1)=\tilde A_n$ and $\tilde A_n\notin W$.
\end{enumerate}
If $A_n\notin W$, we define $I_n$ to be the identity map.  Otherwise we choose a point $\tilde A_n\notin W$ near $A_n$ and draw a smooth curve $I_n$ connecting $A_n$  and $\tilde A_n$. 

Step 2: Extend the  the above homotopies to the subcomplex $\{A_n\}_{n=1}^N\bigcup\{\tau_l\}_{l=1}^L$. To be more precise, we shall construct homotopies $T_l:[0,1]\times [0,1]\rightarrow X$ for $l=1,\cdots,L$ satisfying the following properties:
\begin{enumerate}
	\item  $T_l$ is continuous;
	\item  $T_l(\cdot,0)=\tau_l(\cdot)$;
	\item  $T(0,\cdot)=I_{l_0}(\cdot)$ and $T(1,\cdot)=I_{l_1}(\cdot)$, where $l_0,l_1\in\{1,\cdots,N\}$ such that $\tau_l(0)=A_{l_0}$ and $\tau_l(1)=A_{l_1}$.
\end{enumerate}

\begin{figure}[htbp]
	\centering
	\includegraphics[height=6cm,width=15cm]{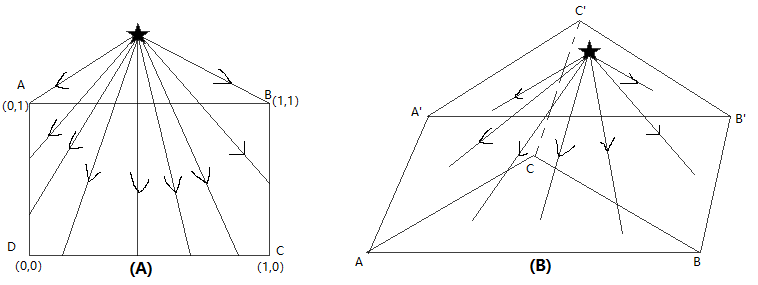}
	\caption{Extension of the homotopies}
	\label{extension}
\end{figure}

Firstly we define the homotopies of $1$-simplexes $\{\tilde{\tau},\tilde\tau_1,\cdots,\tilde{\tau}_m\}$ to be identities. Next we will construct homotopy of $T_l$ for the remaining $\tau_l$ as shown in Figure $(\ref{extension}\text{A})$.  It is easy to see that the projection from the star-shaped point induces a strong deformation retraction $F$ of the unit square to the union of three intervals $\overline{AD}$, $\overline{CD}$ and $\overline{BC}$ in the following sense:
\begin{equation*}
\begin{split}
F:[0,1]\times\big([0,1]\times[0,1]\big)\rightarrow [0,1]\times[0,1],\\
(t,x,y)\mapsto (F_1(t,x,y),F_2(t,x,y)),
\end{split}
\end{equation*}
where $F(0,x,y)=(x,y)$, $F(1,x,y)\subset\overline{AD}\cup\overline{CD}\cup\overline{BC}$ and $F(t,x,y)$ is an identity map when $(x,y)\in\overline{AD}\cup\overline{CD}\cup\overline{BC}.$  Denote by $\widetilde T_l$ the following continuous map defined on $\overline{AD}\cup\overline{CD}\cup\overline{BC}$: 
\begin{equation*}
\begin{split}
&\widetilde T_l:\big(\{0\}\times[0,1]\big)\cup\big([0,1]\times\{0\}\big)\cup\big(\{1\}\times[0,1]\big)\rightarrow X;\\
&\widetilde T_l(x,0)=\tau_l(x),\,\, \widetilde T_l(0,y)=I_{l_0}(y)\,\, {\text{and}}\,\,\widetilde T_l(1,y)=I_{l_1}(y) \,\, {\text{for}}\,\,x,y\in[0,1].
\end{split}
\end{equation*}
Then the desired homotopy $T_l$ can be defined by the following formula:
$$T_l(x,y)=\widetilde T_l(F_1(1,x,y),F_2(1,x,y)).$$

\begin{figure}[htbp]
	\centering
	\includegraphics[height=6cm,width=15cm]{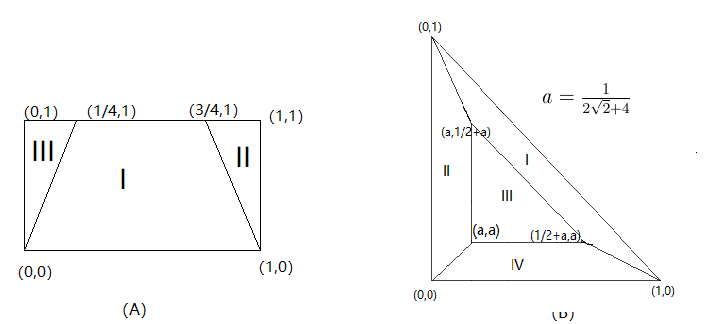}
	\caption{Thickening the boundary}
	\label{thicken}
\end{figure}

Step 3:  Smooth $\{\tau_{l}\}_{l=1}^L$ homotopically with fixed boundary $0$-simplexes $\{A_n\}_{n=1}^N$. To be more precise, we will construct  homotopy $T_l:[0,1]\times [0,1]$ for $l=1,\cdots,L$ satisfying the following properties:
\begin{enumerate}
		\item  $T_l$ is continuous and  $T_l(\cdot,1)$ is a smooth map from $[0,1]$ to $X$;
		\item  $T_l(\cdot,0)=\tau_l(\cdot)$; $T(0,\cdot)=\tau_l(0)$ and $T(1,\cdot)=\tau_l(1)$;
		\item  $T(x,1)=\tau_l(0)$ when $x\approx 0$;  $T(x,1)=\tau_l(1)$  when $x\approx 1$.
\end{enumerate}
Firstly we will thicken the boundary as illustrated by Figure $(\ref{thicken}{\text{A}})$. We define homotopy $T^1_l$ as follows:
\begin{equation*}
\begin{split}
&T^1_l:[0,1]\times[0,1]\rightarrow X;\\
&T^1_l(x,y)=\tau_{l}(0)\,\,{\text{when}}\,\,(x,y)\in{\text{Region III}}\,\,{\text{and}}\,\,T^1_l(x,y)=\tau_{l}(1)\,\,{\text{when}}\,\,(x,y)\in{\text{Region II}};\\
&T^1_l(x,y)=\tau_{l}(\frac{2x-y/2}{2-y})\,\,{\text{when}}\,\,(x,y)\in{\text{Region III}}.
\end{split}
\end{equation*}

Notice that the map $T_l^1(\cdot,1)$ is smooth in a small neighborhood of two boundary points. Then
by Theorem $(10.1.2)$ in \cite{DFN} and Remark $(2)$ therein, we can find a smooth map $\hat{\tau}_l$ homotopic to $T_l^1(\cdot,1)$. Moreover by Remark $(1)$ of Theorem $(10.1.2)$ in \cite {DFN}, $\hat{\tau}_l $ can be chosen so that $\hat{\tau}_l(x)=T_l^1(x,1)$ if $x\approx 0$ or $x\approx 1.$

Step 4: Perturb $\{\tau_l\}_{l=1}^L$ homotopically with fixed boundary points so that $\tau_l$ is disjoint from $W$ for $l=1,\cdots,L.$  First we take a stratification of $W$ as $W=W_1\cup W_2\cup\cdots\cup W_v$, where $W_i$ is smooth for $i=1,\cdots,v$ and ${\text{dim}}W_i<{\text{dim}}W_{i-1}$ for $i=2,\cdots,v.$ According to the proof of Theorem $(10.3.2)$ in \cite{DFN}, we can find a smooth $1$-simplex $\tilde\tau_l^v$ for $\tau_l$ such that
\begin{enumerate}
	\item $\tilde\tau_l^v$ is transversal to $W_v$ ;
	\item  $\tilde\tau_l^v$ is homotopic to $\tau_l$;
	\item  $\tilde\tau_l^v$ coincides with $\tau_l$ in a small neighborhood of boundary points.
\end{enumerate}
By dimension counting, we have that $\tilde\tau_l^v$ is disjoint from $W_v$. Repeating the procedure for $W_{v-1},W_{v-2},\cdots,W_1$, we end up with $\tilde\tau_l^1$ which is disjoint with $W.$ Notice that $\{\tilde{\tau},\tilde\tau_1,\cdots,\tilde{\tau}_m\}$ satisfy the desired property without any perturbation.

Step 5: Extend the obtained homotopies of  $\{\tau_l\}_{l=1}^L$ to $\{\sigma_j\}_{j=1}^J.$ After the above steps we have a homotopy $T_l$ for $\tau_l$ such that 
\begin{enumerate}
	\item  $T_l$ is continuous;
	\item  $T_l(\cdot,0)=\tau_l(\cdot)$, $T(x,1)=T(0,1)$ for $x\approx 0$ and $T(x,1)=T(1,1)$ for $x\approx 1$;
	\item  $T_l(\cdot,1)$ is smooth and disjoint from $W.$
\end{enumerate}

We shall construct homotopy $S_j:\Delta_2\times [0,1]\rightarrow X$ for $j=1,\cdots,J$ satisfying the following properties:
\begin{enumerate}
	\item  $S_j$ is continuous.
	\item  $S_j(\cdot,0)=\sigma_j(\cdot)$.
	\item  Suppose $S_j(0,t,0)=\pm\tau_{j_{01}}(t)$,
	$S_j(t,0,0)=\pm\tau_{j_{02}}(t)$ and
	$S_j(t,1-t,0)=\pm\tau_{j_{12}}(t)$ (the signs depend on the orientation) for $j_{01},j_{02},j_{12}\in\{1,\cdots,J\}$ and $t\in[0,1]$. Then
	 $S_j(0,t,\cdot)=\pm T_{j_{01}}(t,\cdot)$,
	$S_j(t,0,\cdot)=\pm T_{j_{02}}(t,\cdot)$ and
	$S_j(t,1-t,\cdot)=\pm T_{j_{12}}(t,\cdot)$ for $t\in[0,1]$ with the compatible signs. 
\end{enumerate}

Similar to Step 2, we first construct a strong deformation retract as illustrated in Figure $(\ref{extension}{\text{B}})$. It is clear that the projection from the star-shaped point induces a strong deformation retraction $F$ of the triangle based prism to
 to the union of the triangle $\overline{ABC}$ and three parallelograms $\overline{ACC^{\prime}A^{\prime}}$, $\overline{CBB^{\prime}C^{\prime}}$ and $\overline{ACC^{\prime}A^{\prime}}$ in the following sense:
\begin{equation*}
\begin{split}
F:[0,1]\times\big(\Delta_2\times[0,1]\big)&\rightarrow \Delta_2\times[0,1],\\
(t,x,y,s)&\mapsto (F_1(t,x,y,s),F_2(t,x,y,s),F_3(t,x,y,s)),
\end{split}
\end{equation*}
where $F(0,x,y,t)=(x,y,t)$, $F(1,x,y,t)\subset\overline{AD}\cup\overline{CD}\cup\overline{BC}$ and $F(t,x,y,s)$ is an identity map when $(x,y,s)\in\overline{ABC}\cup\overline{ACC^{\prime}A^{\prime}}\cup\overline{CBB^{\prime}C^{\prime}}\cup\overline{ACC^{\prime}A^{\prime}}.$ 

We construct $S_j$ for $\sigma_j$ as follows. Denote by $\widetilde S_j$ the following continuous map defined on $\overline{AD}\cup\overline{CD}\cup\overline{BC}$: 
\begin{equation*}
\begin{split}
&\widetilde S_j:\overline{ABC}\cup\overline{ACC^{\prime}A^{\prime}}\cup\overline{CBB^{\prime}C^{\prime}}\cup\overline{ACC^{\prime}A^{\prime}}\rightarrow X;\\
&\widetilde S_j|_{\overline{ABC}}=\sigma_j,\,\, \widetilde S_j|_{ACC^{\prime}A^{\prime}}=\pm T_{j_{01}},\,\,\widetilde S_j|_{CBB^{\prime}C^{\prime}}=\pm T_{j_{02}}\,\, {\text{and}}\,\,\widetilde S_j|_{ACC^{\prime}A^{\prime}}=\pm T_{j_{12}}.
\end{split}
\end{equation*}
Then the desired homotopy $S_j$ is given by
$$S_j(x,y,s)=\widetilde S_j(F_1(1,x,y,s),F_2(1,x,y,s),F_3(1,x,y,s)).$$

Step 6:  Smooth $\{\sigma_{j}\}_{j=1}^J$ homotopically with fixed boundary $1$-simplexes $\{\tau_l\}_{l=1}^L$. To be more precise, we will construct  homotopy $S_j:\Delta_2\times [0,1]$ for $j=1,\cdots,J$ satisfying the following properties:
\begin{enumerate}
	\item  $S_j$ is continuous and $S_j(\cdot,1)$ is a smooth map from $\Delta_2$ to $X$.
	\item  $S_j(\cdot,0)=\sigma_j(\cdot)$. 
	\item  $S_j(x,y,\cdot)=\sigma_j(x,y)$  for $(x,y)\in\big(\{0\}\times[0,1]\big)\cup\big([0,1]\times\{0\}\big)\cup\{(t,1-t)|0\leq t\leq 1\}$.
\end{enumerate}
We proceed in a similar way to Step 3. Firstly we thicken the boundary as illustrated by Figure $(\ref{thicken}{\text{B}})$. We can define a homotopy $S^1_j:\Delta_2\times[0,1]\rightarrow X$ such that $S^1_j(\cdot,\cdot,1)$ is defined as follows:
\begin{equation*}
\begin{cases}
S^1_j(x,y,1)=\sigma_j(\frac{\frac{1}{\sqrt2+2}y+(1-\frac{1}{\sqrt2+2})x-\frac{1}{\sqrt2+2}}{x+y-\frac{2}{\sqrt2+2}},\frac{\frac{1}{\sqrt2+2}x+(1-\frac{1}{\sqrt2+2})y-\frac{1}{\sqrt2+2}}{x+y-\frac{2}{\sqrt2+2}})\,\,{\text{when}}\,\,(x,y)\in{\text{Region I}};\\
S^1_j(x,y,1)=\sigma_j(2x-\frac{1}{\sqrt2+2},2y-\frac{1}{\sqrt2+2})\,\,{\text{when}}\,\,(x,y)\in{\text{Region II}}.\\
S^1_j(x,y,1)=\sigma_j(0,\frac{1}{\sqrt2+2}\frac{x-y}{x-\frac{1}{\sqrt2+2}})\,\,{\text{when}}\,\,(x,y)\in{\text{Region III}};\\
S^1_j(x,y,1)=\sigma_j(\frac{1}{\sqrt2+2}\frac{y-x}{y-\frac{1}{\sqrt2+2}},0)\,\,{\text{when}}\,\,(x,y)\in{\text{Region IV}}.\\
\end{cases}
\end{equation*}
Notice that after Step 1-5 $\tau_j$ is smooth and constant near the  boundary $\tau_j(0)$ and $\tau_j(1)$. Then $S^1_j(\cdot,\cdot,1)$ is smooth in a small neighborhood of the boundary.  By Theorem $(10.1.2)$ in \cite{DFN} and   Remark $(2)$ of it, we can find a smooth map $\hat{\sigma}_j$ homotopic to $S^1_j(\cdot,\cdot,1)$. Moreover by Remark $(1)$ of Theorem $(10.1.2)$ in \cite{DFN},  we can assume  $\hat{\sigma}_j$ coincides with  $\tilde{\sigma}_j$ in a small neighborhood of the boundary of $\Delta_2$. Hence we get the desired homotopy of $\sigma_j$.

Step 7: Perturb $\{\sigma_j\}_{j=1}^J$ so that $\sigma_j$ is disjoint from the singularities of $W$ and intersects the smooth part of $W$ transversally.  First we take a stratification of $W$ as $W=W_1\cup W_2\cup\cdots\cup W_v$, where $W_i$ is smooth for $i=1,\cdots,v$ and  ${\text{dim}}W_i<{\text{dim}}W_{i-1}$ for $i=2,\cdots,v.$ 

Next we will find a desired perturbation for $\sigma_l$. Notice that a certain neighborhood of the boundary of $\sigma_j$ is transversal to $W_v$. Then following the proof of Theorem $(10.3.2)$ in \cite{DFN}, we can find a smooth $2$-simplex $\tilde\sigma_j^v$ for $\sigma_j$ such that $\tilde\sigma_j^v$ is transversal to $W_v$ and homotopic to $\tau_l$ with boundary points fixed. Repeating the procedure for $W_{v-1},W_{v-2},\cdots,W_1$, we derive the corresponding transversal maps $\tilde\sigma_j^{v-1},\tilde\sigma_j^{v-2},\cdots$ and $\tilde\sigma_j^1$ accordingly.  By dimension counting, we conclude that the $\tilde\sigma_j^1$ is disjoint from $\bigcup_{i=2}^v W_{i}$ and intersect $W_1$ transversally in the interior.

It is easy to verify that $\tilde\sigma_j^1$ satisfy all the required property in the lemma.  Hence we complete the proof. $\endpf$

\bigskip
\noindent H. Fang, Department of Mathematics, University of Wisconsin-Madison,
 Madison, WI 53706, USA.  (hfang35$ @$wisc.edu)

\end{document}